\documentclass[11pt]{article}

\usepackage[T1]{fontenc}
\usepackage[utf8]{inputenc}
\usepackage{amsmath,amssymb,epsfig,bbm}
\usepackage{stmaryrd}
\usepackage[colorlinks]{hyperref}
\usepackage{comment}
\usepackage{color}
\usepackage{textcomp}

\usepackage[acronym,toc]{glossaries}
\usepackage{glossary-longragged}
\usepackage{glossary-superragged}
\newglossary[slg]{symbols}{sym}{sbl}{List of Symbols}

\usepackage[textsize=small]{todonotes}
\usepackage{varwidth}

\usepackage{url}

\usepackage{enumitem}
\setlist[itemize]{labelindent=*, leftmargin=.5 truecm,nosep}

\usepackage[nottoc]{tocbibind}

\usepackage{fancyhdr}
\usepackage[us,12hr]{datetime} 
\fancypagestyle{plain}{
\fancyhf{}
\rfoot{ {\ddmmyyyydate\today}}
\lfoot{\thepage}
}

\usepackage{amsthm}

\theoremstyle{plain}
\newtheorem{defi}{Definition}[section]
\newtheorem{prop}[defi]{Proposition}
\newtheorem{theo}[defi]{Theorem}
\newtheorem{conj}[defi]{Conjecture}
\newtheorem{lemm}[defi]{Lemma}
\newtheorem{coro}[defi]{Corollary}

\theoremstyle{definition}
\newtheorem{rema}[defi]{Remark}
\newtheorem{exem}[defi]{Example}
\newtheorem{exems}[defi]{Examples}

\newcommand{\bdefi}{\begin{defi}}
\newcommand{\edefi}{\end{defi}}
\newcommand{\bprop}{\begin{prop}}
\newcommand{\eprop}{\end{prop}}
\newcommand{\btheo}{\begin{theo}}
\newcommand{\etheo}{\end{theo}}
\newcommand{\blemm}{\begin{lemm}}
\newcommand{\brema}{\begin{rema}}
\newcommand{\erema}{\end{rema}}
\newcommand{\bexer}{\begin{exem}}
\newcommand{\eexer}{\end{exem}}
\newcommand{\bexem}{\begin{exem}}
\newcommand{\eexem}{\end{exem}}
\newcommand{\bexems}{\begin{exems}}
\newcommand{\eexems}{\end{exems}}
\newcommand{\bconj}{\begin{conj}}
\newcommand{\econj}{\end{conj}}
\newcommand{\elemm}{\end{lemm}}
\newcommand{\bcoro}{\begin{coro}}
\newcommand{\ecoro}{\end{coro}}
\newcommand{\dem}{\noindent{\bf Proof. }}

\usepackage{mathrsfs}
\renewcommand\mathcal{\mathscr}

\newcommand{\E}{{\cal E}}
\newcommand{\F}{{\cal F}}
\newcommand{\G}{{\cal G}}
\renewcommand{\H}{{\cal H}}

\newcommand{\OOO}{{\cal O}}

\newcommand{\maths}[1]{{\mathbb #1}}  

\newcommand{\FF}{\maths{F}}

\newcommand{\NN}{\maths{N}}

\newcommand{\PP}{\maths{P}}
\newcommand{\QQ}{\maths{Q}}
\newcommand{\RR}{\maths{R}}

\newcommand{\TT}{\maths{T}}

\newcommand{\ZZ}{\maths{Z}}

\newcommand{\mmm}{{\mathfrak m}}

\newcommand{\bs}{\backslash}
\newcommand{\ga}{\gamma}
\newcommand{\Ga}{\Gamma}
\newcommand{\ov}[1]{{\overline{#1}}} 
\newcommand{\ra}{\rightarrow}
\newcommand{\wt}[1]{{\widetilde{#1}}}
\newcommand{\wh}[1]{{\widehat{#1}}}

\newcommand{\Aut}{\operatorname{Aut}}
\newcommand{\Ax}{\operatorname{Ax}}

\newcommand{\card}{{\operatorname{Card}}}

\newcommand{\cqfd}{\hfill$\Box$}

\newcommand{\dbs}{\backslash\!\!\backslash}

\newcommand{\diam}{{\operatorname{diam}}}

\newcommand{\gengeod}
{\operatorname{\widecheck{\G\,}\!\!}}

\newcommand{\length}{\operatorname{length}}

\newcommand{\lk}{\operatorname{lk}}

\newcommand{\QI}{K^{(2)}}

\newcommand{\Sp}{\operatorname{Sp}}

\newcommand{\GL}{\operatorname{GL}}
\newcommand{\PGL}{\operatorname{PGL}}

\newcommand{\PSL}{\operatorname{PSL}}

\usepackage{interval}
\newcommand{\cli}{\interval} 
\newcommand{\oi}{\interval[open]} 
\newcommand{\cloi}{\interval[open right]} 

\usepackage{mathtools}
\makeatletter
\DeclareRobustCommand\widecheck[1]{{\mathpalette\@widecheck{#1}}}
\def\@widecheck#1#2{%
    \setbox\z@\hbox{\m@th$#1#2$}%
    \setbox\tw@\hbox{\m@th$#1%
       \widehat{%
          \vrule\@width\z@\@height\ht\z@
          \vrule\@height\z@\@width\wd\z@}$}%
    \dp\tw@-\ht\z@
    \@tempdima\ht\z@ \advance\@tempdima2\ht\tw@ \divide\@tempdima\thr@@
    \setbox\tw@\hbox{%
       \raise\@tempdima\hbox{\scalebox{1}[-1]{\lower\@tempdima\box
\tw@}}}%
    {\ooalign{\box\tw@ \cr \box\z@}}}
\makeatother

\newcounter{fig}

\title{On the nonarchimedean quadratic Lagrange spectra}
\author{Jouni Parkkonen \and Fr\'ed\'eric Paulin} 

\begin{document}
\bibliographystyle{../alphanum}
\maketitle
\begin{abstract}
We study Diophantine approximation in completions of functions fields
over finite fields, and in particular in fields of formal Laurent
series over finite fields. We introduce a Lagrange spectrum for the
approximation by orbits of quadratic irrationals under the modular
group. We give nonarchimedean analogs of various well known results
in the real case: the closedness and boundedness of the Lagrange
spectrum, the existence of a Hall ray, as well as computations of
various Hurwitz constants. We use geometric methods of group actions
on Bruhat-Tits trees.
\end{abstract}

\section{Introduction}
\label{sec:intro}

Diophantine approximation in local fields of positive characteristic
has been studied for many years, both from the classical viewpoint
(see the works of de Mathan, Lasjaunias, Osgood, W.~Schmidt, Thakur,
Voloch for instance) and from the point of view of arithmetic geometry
(see the works of Mahler, Manin, and many others), as well from an
ergodic theory viewpoint (see for instance \cite{BerNak00}). In this
paper, we are interested in the approximation by quadratic
irrationals: we define and study the {\em quadratic} Lagrange
spectra in completions of function fields over finite fields with
respect to the absolute values defined by discrete valuations.  In
this introduction, we concentrate on the special case of the field of
rational fractions and its valuation at infinity. We refer to Section
\ref{sec:hall} for the definitions and results in the general case,
allowing congruence considerations.

\medskip
Let $\FF_q$ be a finite field of order a positive power $q$ of a
positive prime.  Let $R=\FF_q[Y]$, $K=\FF_q(Y)$ and $\wh{K}=
\FF_q((Y^{-1}))$ be respectively the ring of polynomials in one
variable $Y$ over $\FF_q$, the field of rational functions in $Y$ over
$\FF_q$ and the field of formal Laurent series in $Y^{-1}$ over
$\FF_q$.  Then $\wh{K}$ is a nonarchimedean local field, the
completion of $K$ with respect to its place at infinity, that is, the
absolute value $|\frac{P}{Q}|=q^{\deg P-\deg Q}$ for all $P,Q\in
R-\{0\}$.  Let
$$
\QI=\{f\in \wh{K}\;:\;[K(f):K]=2\}
$$ 
be the set of quadratic irrationals over $K$ in $\wh{K}$.  Given $f\in
\wh{K}-K$, it is well known that $f\in \QI$ if and only if the
continued fraction expansion\footnote{See Section \ref{sec:comput} for
  a definition.} of $f$ is eventually periodic. The projective action
of $\Ga=\PGL_2(R)$ on $\PP_1(\wh K)=\wh K\cup \{\infty\}$ preserves
$\QI$, keeping the periodic part of the continued fraction expansions
unchanged (up to cyclic permutation and invertible elements). We refer
for instance to \cite{Lasjaunias00, Schmidt00, Paulin02} for
background on the above notions.

Now let us fix $\alpha\in \QI$.  We denote by $\alpha^\sigma\in \QI$
the Galois conjugate of $\alpha$ over $K$.  The {\it complexity}
$h(\alpha)=\frac{1}{|\alpha-\alpha^\sigma|}$ of $\alpha$ was
introduced in \cite{HerPau10} and developped in \cite[\S
  17.2]{BroParPau18}. It plays the role of the (naive) height of a
rational number in Diophantine approximation by rationals, and is an
appropriate complexity when studying the approximation by elements in
the orbit under the modular group of a given quadratic
irrational.\footnote{see \cite{ParPau11MZ} for a comparison between
  the naive height of the algebraic number $\alpha$ and the above
  complexity of $\alpha$ in the Archimedean case.} We refer to the
above references for motivations and results, in particular to
\cite[Thm.~1.6]{HerPau10} for a Khintchine type result and to \cite[\S
  17.2]{BroParPau18} for an equidistribution result of the orbit of
$\alpha$ under $\PGL_2(R)$.

Let
$$
\Theta_\alpha=\PGL_2(R)\cdot\{\alpha,\alpha^\sigma\}
$$
be the union of the orbits of $\alpha$ and $\alpha^\sigma$ under the
projective action of $\PGL_2(R)$. Given $x\in\wh K - (K\cup
\Theta_\alpha)$, we define the {\it quadratic approximation constant}
of $x$ by
$$
c_\alpha(x)=\liminf_{\beta\in\Theta_\alpha,\;|\beta-\beta^\sigma|\ra 0}
\;\frac{|x-\beta|}{|\beta-\beta^\sigma|}\;.
$$
We define the {\it quadratic Lagrange spectrum} of $\alpha$ as
$$
\Sp(\alpha)=\{c_\alpha(x)\;:\;x\in\wh K - (K\cup \Theta_\alpha)\}\;.
$$ 
Note that $\Sp(\alpha)\subset q^\ZZ\cup \{0,+\infty\}$. It follows
from \cite[Thm.~1.6]{HerPau10} that if $m_{\wh K}$ is a Haar measure
on the locally compact additive group of $\wh K$, then for $m_{\wh
  K}$-almost every $x\in\wh K$, we have $c_\alpha(x)=0$. Hence in
particular, $0\in \Sp(\alpha)$ and the quadratic Lagrange spectrum is
therefore closed. In Section \ref{sec:hall}, we prove that it is
bounded, and we can thus define the {\it (quadratic) Hurwitz
  constant} of $\alpha$ as $\max\Sp(\alpha)\;\in q^\ZZ$.

The following theorems, giving nonarchimedean analogs of the results
of Lin, Bugeaud and Pejkovic \cite{Lin18,Bugeaud14,Pejkovic16}, say
that the quadratic Lagrange spectrum of $\alpha$ is a closed bounded
subset of $q^\ZZ\cup \{0\}$ which contains an initial interval, and
computes various Hurwitz constants.

\btheo\label{theo:hurwitzhallintro} Let $\alpha$ be a quadratic irrational
over $K$ in $\wh{K}$.

(1) (Upper bound) Its quadratic Hurwitz constant satisfies
$\max\Sp(\alpha)\leq q^{-2}$.

(2) (Hall ray) There exists $m_\alpha\in\NN$ such that for all
$n\in\NN$ with $n\geq m_\alpha$, we have $q^{-n}\in\Sp(\alpha)$.
\etheo

In Section \ref{sec:hall}, we even prove that Assertion (2) of this
theorem is valid when $K$ is any function field over $\FF_q$, $\wh K$
is the completion of $K$ at any place of $K$, and $R$ is the
corresponding affine function ring.

\btheo\label{theo:maxhurwitzintro} The Hurwitz constant of any
quadratic irrational over $K$ in $\wh{K}$, whose continued fraction
expansion is eventually $k$-periodic with $k\leq q-1$, is equal to
$q^{-2}$.  \etheo

There are examples of quadratic irrationals for which the quadratic
Lagrange spectrum coincides with the maximal Hall ray. The following
theorem gives a special case, see Theorem \ref{theo:goldenspec} for a
more general result.

\btheo\label{theo:goldenratiointro} If 
$$
\varphi =Y + \cfrac{1}{Y+\cfrac{1}{Y+\cdots}}\,,
$$ 
then $\Sp(\varphi)=\{0\}\cup\{q^{-n-2}\;:\;n\in\NN\}$. \etheo 

In Proposition \ref{prop:gaps}, we give a class of quadratic
irrationals whose quadratic Lagrange spectrum does not coincide with
its maximal Hall ray, in other words, who have {\em gaps} in their
spectrum.

\medskip After the first version of this paper was posted on ArXiv,
Yann Bugeaud \cite{Bugeaud18} has given a completely different proof
of the above results (except the generalisation to function fields),
and proved several new theorems giving a more precise description of
these spectra. In particular, he proved that all approximation
constants for a given quadratic irrational are attained on the other
quadratic irrationals, that for every $m\geq 2$ there exists $\beta\in
\QI$ such that $\max\Sp(\beta)=q^{-m}$, and that for all
$\ell\in\NN$, there exists $\beta\in \QI$ such that $\Sp(\beta)$
contains exactly $\ell$ gaps.

\medskip
In order to explain the origin of our results, recall that for
$x\in\RR-\QQ$, the {\it approximation constant} of $x$ by rational
numbers is
$$
c(x)=
\;\liminf_{p,q\in\ZZ,\; q\ra+\infty}\;\; q^2\Big|x-\frac pq\Big|,
$$ and that the {\it Lagrange spectrum} is $\Sp_\QQ= \{c(x)\;:\; x\in
\RR-\QQ\}$. Numerous properties of the Lagrange spectrum are known,
see for instance \cite{CusFla89}. In particular, $\Sp_\QQ$ is bounded
and closed, has maximum $\frac{1}{\sqrt{5}}$, and contains a maximal
interval $\interval0\mu$ with $0<\mu<\frac{1}{\sqrt{5}}$ called a {\it
  Hall ray}.  Khinchin \cite{Khinchin64} proved that almost every real
number is badly approximable by rational numbers, so that the
approximation constant vanishes almost surely.  Many of these results
have been generalised to the Diophantine approximation of complex
numbers, Hamiltonian quaternions and for the Heisenberg group, see for
example
\cite{Poitou53,Schmidt69,Schmidt75,ParPau07,ParPau09,ParPau10GT}.

Let $\alpha_0$ be a fixed real quadratic irrational number over
$\QQ$. For every such number $\alpha$, let $\alpha^\sigma$ be its
Galois conjugate. Let $\E_{\alpha_0}=\PSL_2(\ZZ)\cdot \{\alpha_0,
\alpha_0^\sigma\}$ be its (countable, dense in $\RR$) orbit for the
action by homographies and anti-homographies of ${\rm PSL}_2(\ZZ)$ on
$\RR\cup\{\infty\}$.  For every $x\in \RR-(\QQ\cup \E_{\alpha_0})$,
the {\it approximation constant} of $x$ by elements of $\E_{\alpha_0}$
was defined in \cite{ParPau11MZ} by
$$
c_{\alpha_0}(x)=
\;\liminf_{\alpha\in\E_{\alpha_0}\;:\;|\alpha-\alpha^\sigma|\ra 0}\;\;
2\,\frac{|x-\alpha|}{|\alpha-\alpha^\sigma|}\;,
$$
the {\it quadratic Lagrange spectrum} (or {\it approximation
  spectrum}) of $\alpha_0$ by
$$
\Sp(\alpha_0) =\{c_{\alpha_0}(x)\;:\;x\in \RR-(\QQ\cup
\E_{\alpha_0})\}\;,
$$ 
and the {\it Hurwitz constant} of $\alpha_0$ by $\sup\Sp(\alpha_0)$.
We proved that the quadratic Lagrange spectrum of $\alpha_0$ is
bounded and closed, and that an analog of Khinchin's theorem holds. We
generalised the definitions and the above results to the approximation
of complex numbers and elements of the Heisenberg group.  In the
latter cases, we also proved the existence of a Hall ray in the
spectrum.

In the real case, the existence of a Hall ray in $\Sp(\alpha_0)$ is
due to Lin \cite{Lin18}.  Bugeaud \cite{Bugeaud14} proved that the
Hurwitz constant of the Golden Ratio $\phi$ is equal to $\frac 3{\sqrt
  5}-1$, and his conjecture that the Hurwitz constant of any real
quadratic irrational is at most $\frac 3{\sqrt 5}-1$ was confirmed by
Pejkovic \cite{Pejkovic16}. The Hurwitz constant is known explicitly
in many $2$-periodic continued fraction expansion cases, see
\cite{Pejkovic16,Lin18}.


\medskip\noindent{\small {\it Acknowledgements: } This work was
  supported by the French-Finnish CNRS grant PICS 
  \textnumero\,6950. We thank a lot Yann Bugeaud for his comments on
  the first version of this paper, which allowed us to remove the
  unnecessary assumption of odd characteristic, and gave a negative
  solution to a conjecture we proposed on a general formula for the
  Hurwitz constants. }

\section{Background on function fields and Bruhat-Tits trees} 
\label{sec:background}

In this section, we recall the basic notations and properties of
function fields $K$ over $\FF_q$ and their valuations $v$, the
associated Bruhat-Tits trees $\TT_v$ and modular groups $\Ga_v$ acting
on $\TT_v$.  We refer to \cite{Goss96, Rosen02, Serre83} for
definitions, proofs and further information, see also \cite[Ch.~14 and
  15]{BroParPau18}.

Let $\FF_q$ be a finite field of order $q$ with $q$ a positive power
of a positive prime.  

\subsection{Function fields}

Let $K$ be a function field over $\FF_q$ and let $v:K^\times \ra\ZZ$
be a (normalised discrete) valuation of $K$. Let $R_v$ be the affine
function ring associated with $(K,v)$. Let $|\cdot|_v$ be the absolute
value on $K$ corresponding to $v$ and let $K_v$ be the completion of
$K$ with respect to $|\cdot|_v$. We again denote by $v$ and
$|\cdot|_v$ the extensions of $v$ and $|\cdot|_v$ to $K_v$. Let
$$
\OOO_v=\{x\in K_v\;:\;v(x)\ge 0\}
$$ 
be the {\em valuation ring} of $K_v$.  Its unique maximal ideal is
$$
\mathfrak m_v=\{x\in K_v\;:\;v(x)> 0\}\,.
$$ 
We denote the cardinality of the {\em residual field} $k_v=
\OOO_v/\mathfrak m_v$ by $q_v$, so that $|f|_v=q_v^{-v(f)}$ for all
$f\in K_v^\times$.

\bexem 
Let $K=\FF_q(Y)$ be the field of rational functions in one
variable $Y$ with coefficients in $\FF_q$, let $\FF_q[Y]$ be the ring
of polynomials in one variable $Y$ with coefficients in $\FF_q$, and
let $v_\infty: K^\times \ra \ZZ$ be the {\em valuation at
  infinity}\index{valuation!at infinity} of $K$, defined on every
$P/Q\in K$ with $P,Q\in\FF_q[Y]-\{0\}$ by
$$
v_\infty(P/Q)=\deg Q-\deg P\;.
$$ Then $R_{v_\infty}=\FF_q[Y]$ and the absolute value $|\cdot|_{v_\infty}$ on
$K$ associated with $v_\infty$ is the absolute value defined in the
introduction.  The completion $K_{v_\infty}$ of $K$ for
$|\cdot|_{v_\infty}$ is the field $\FF_q((Y^{-1}))$ of formal Laurent
series in one variable $Y^{-1}$ with coefficients in $\FF_q$, denoted
by $\wh K$ in the introduction. The elements $x$ in $\FF_q((Y^{-1}))$
are of the form
$$
x= \sum_{i\in\ZZ} x_i\,Y^{-i}
$$ where $x_i\in\FF_q$ for all $i\in\ZZ$ and $x_i=0$ for
$i$ small enough. The valuation at infinity of $\FF_q((Y^{-1}))$
extending the valuation at infinity of $\FF_q(Y)$ is
$$
\forall\;x\in \FF_q((Y^{-1}))^\times,\;\;\;
v_\infty(x)=\sup\{i\in\ZZ\;:\;\forall\;j<i,\;\;\;x_j=0\}\,.
$$ 
The valuation ring $\OOO_{v_\infty}$ of $v_\infty$ is the ring
$\FF_q[[Y^{-1}]]$ of formal power series in one variable $Y^{-1}$ with
coefficients in $\FF_q$. The residual field $k_{v_\infty}=
\OOO_{v_\infty}/ \mmm_{v_\infty}$ of $v_\infty$ is
$\FF_q$, which has order $q_{v_\infty}=q$.  
\eexem

We identify the projective line $\PP_1({K}_v)$ with ${K}_v \cup
\{\infty\}$ using the map ${K}_v(x,y) \mapsto \frac{x}{y}$, so that
$$
\infty=[1:0]\;.
$$ 
The projective action of $\PGL_2({K}_v)$ on $\PP^1({K}_v)$ is the
action by homographies on ${K}_v\cup\{\infty\}$, given by
$(g,z)\mapsto g\cdot z =\frac{a\,z+b}{c\,z+d}$ if $g=\begin{bmatrix} a
& b \\ c & d\end{bmatrix}\in\PGL_2({K}_v)$. As usual, we define
$g\cdot\infty=\frac ac$ and $g\cdot(-\frac dc)=\infty$.

\subsection{Bruhat-Tits trees}

An {\it $\OOO_{v}$-lattice}\index{lattice@$\OOO_{v}$-lattice}
$\Lambda$ in the ${K}_v$-vector space ${K}_v\times {K}_v$ is a rank
$2$ free $\OOO_{v}$-submodule of ${K}_v\times {K}_v$, generating
${K}_v\times {K}_v$ as a vector space. The Bruhat-Tits tree $\TT_v$ of
$(\PGL_2,{K}_v)$ is the graph whose set of vertices $V\TT_v$ is the
set of homothety classes (under $({K}_v)^\times$) $[\Lambda]$ of
$\OOO_{v}$-lattices $\Lambda$ in ${K}_v\times {K}_v$, and whose set of
edges $E\TT_v$ is the set of pairs $(x,x')$ of vertices such that
there exist representatives $\Lambda$ of $x$ and $\Lambda'$ of $x'$
for which $\Lambda \subset \Lambda'$ and $\Lambda'/\Lambda$ is
isomorphic to $\OOO_{v}/\mmm_{v}$.  The graph $\TT_v$ is a
regular tree of degree $|\PP_1(k_{v})|= q_v+1$.  The {\it standard
  base point}\index{standard base point} $*_v$ of $\TT_v$ is the
homothety class $[\OOO_v\times\OOO_v]$ of the $\OOO_{v}$-lattice
$\OOO_{v} \times \OOO_{v}$ generated by the canonical basis of ${K}_v
\times {K}_v$.
The {\it link}
$$
\lk(*_v)=\{y\in V\TT_v\;:\;d(y,*_v)=1\}
$$ 
of $*_v$ in $\TT_v$ identifies with the projective line $\PP_1(k_v)$.

The left linear action of $\GL_2({K}_v)$ on ${K}_v\times {K}_v$
induces a faithful, vertex-transitive left action of $\PGL_2 ({K}_v)$
by automorphisms on $\TT_v$.  The stabiliser of $*_v$ in $\PGL_2
({K}_v)$ is $\PGL_2 (\OOO_v)$, which acts projectively on $\lk(*_v)=
\PP_1(k_v)$ by reduction modulo $v$, and in particular $\PGL_2(k_v)$
acts simply transitively on triples of pairwise distinct points on
$\lk(*_v)$. We identify the boundary at infinity $\partial_\infty
\TT_v$ of $\TT_v$ with $\PP_1({K}_v)$ by the unique homeomorphism from
$\partial_\infty \TT_v$ to $\PP_1({K}_v)$ such that the (continuous)
extension to $\partial_\infty \TT_v$ of the isometric action of
$\PGL_2({K}_v)$ on $\TT_v$ corresponds to the projective action of
$\PGL_2({K}_v)$ on $\PP_1({K}_v)$.

Let $\Ga_v=\PGL_2(R_v)$.  The group $\Ga_v$ is a lattice in the
locally compact group $\PGL_2({K}_v)$, called the {\it modular group}
at $v$ of $K$.  The quotient graph $\Ga\bs\TT_v$ is called the {\it
  modular graph} of $K$, and the quotient graph of groups
$\Ga\dbs\TT_v$ is called the {\it modular graph of groups} at $v$ of
$K$.  We refer to \cite{Serre83} for background information on these
objects, and for instance to \cite{Paulin02} for a geometric treatment
when $K=\FF_q(Y)$ and $v= v_\infty$.

Recall that the {\em open horoballs} centred at $\xi\in
\partial_{\infty} \TT_v$ are the subsets of the geometric realisation
$|\TT_v|$ of $\TT_v$ defined by
$$
\H(\rho_\xi)=\{y\in|\TT_v|\;:\; \lim_{t\to+\infty}\big(
t-d(\rho_{\xi}(t),y) \big)> 0\}
$$ 
where $\rho_{\xi}$ is a geodesic ray converging to $\xi$. The boundary
of $\H(\rho_\xi)$ is the {\em horosphere} 
$$
\partial\H(\rho_\xi)=\{y\in|\TT_v| \;:\;
\lim_{t\to+\infty} \big(t-d(\rho_{\xi}(t),y)\big)= 0\}\,.
$$ 
We refer to \cite{BriHae99} for background on these notions.  The {\em
  height} in $\H(\rho_\xi)$ of a point $x\in|\TT_v|$ is
$\lim_{t\to+\infty} \big(t-d(\rho_{\xi}(t),x)\big)$. It is positive if
and only if $x$ belongs to $\H(\rho_\xi)$.  We denote by $\H_\infty$
the unique horoball centred at $\infty\in\partial_\infty \TT_v$ whose
associated horosphere passes through $*_v$.

Let $\Ga$ be a finite index subgroup of $\Ga_v$. By for instance
\cite{Serre83, Paulin02}, there exists a $\Ga$-equivariant family of
pairwise disjoint open horoballs $(\H_\xi)_{\xi\in\PP_1(K)}$ in
$\TT_v$ with $\H_\xi$ centered at $\xi$ and the stabiliser $\Ga_\xi$
of $\xi$ in $\Ga$ acting transitively on the boundary of $\H_\xi$ for
every $\xi\in\PP_1(K)\subset \partial_\infty \TT_v$, so that the
quotient by $\Ga$ of
$$
\wt E_\Ga=\TT_v-\bigcup_{\xi\in\PP_1(K)}\H_\xi
$$ 
is a finite connected graph, denoted by $E_\Ga$.  The {\it set of
  cusps} $\Ga\bs\PP_1(K)$ is finite.  For every representative $\xi$
of a cusp in $\Ga\bs\PP_1(K)$, the injective image by the canonical
projection $\TT_v\ra \Ga\bs \TT_v$ of any geodesic ray starting from a
point on the boundary of $\H_\xi$ with point at infinity $\xi$ is
called a {\em cuspidal ray}.  The quotient graph $\Ga\bs \TT_\infty$
is the union of the finite subgraph $E_\Ga$ and the finite collection
of cuspidal rays that are glued to $E_\Ga$ at their origin.

\bexem\label{exem:2.2} (See for instance \cite{BasLub01}.) 
Let $K=\FF_q(Y)$ and $v= v_\infty$. Then
$\Ga_{v_\infty}=\PGL_2(\FF_q[Y])$ and the quotient graph of groups
$\Ga_{v_\infty}\dbs\TT_{v_\infty}$ is the following {\it modular ray}

\begin{center}
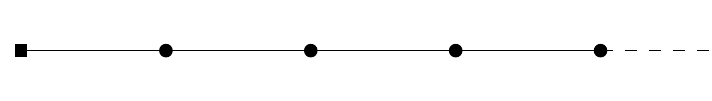
\end{center}

\noindent 
where $\Ga_{-1}=\PGL_2(\FF_q)$, $\Ga'_0=\Ga_0\cap\Ga_{-1}$ and, for
every $n\in\NN$, 
$$
\Ga_n=\bigg\{\begin{bmatrix}a&b\\0&d\end{bmatrix}
\in\PGL_2(\FF_q[Y])
\;:\; a,d\in \FF_q^\times, b\in \FF_q[Y], \deg b\leq n+1\bigg\}\;.
$$ The {\rm canonical} $\Ga_{v_\infty}$-equivariant family
$(\H_\xi)_{\xi\in\PP_1(K)}$ of pairwise disjoint maximal open
horo\-balls in $\TT_{v_\infty}$, with $\H_\xi$ centered at $\xi$ for
every $\xi\in\PP_1(K)$, consists of the connected components of
$|\TT_{v_\infty}| -\Ga_{v_\infty} *_{v_\infty}$. The graph $\wt E_\Ga$
is reduced to the orbit of the base point $*_{v_\infty}$, and $E_\Ga$
is reduced to one vertex, the origin of the modular ray (see the above
picture). In particular, the diameter of $E_\Ga$ is zero.

All geodesic lines in $\TT_{v_\infty}$ pass through the
$\Ga_{v_\infty}$-orbit of $*_{v_\infty}$. Indeed, no geodesic is completely
contained in a horoball and since $\bigcup_{\xi\in\PP_1(K)}\partial
\H_\xi = \Ga_{v_\infty} *_{v_\infty}$, the only way a geodesic line
exits a horoball of the canonical family $(\H_\xi)_{\xi\in\PP_1(K)}$
is through the orbit of $*_{v_\infty}$.  The intersection of a
geodesic line with the orbit $\Ga_{v_\infty} *_{v_\infty}$ is hence
finite if and only if its endpoints are both in $\PP_1(K)$.  
\eexem

\medskip 
We end this section with the following lemma, which is an effective
version of a special case of \cite[Prop.~2.6]{HerPau10}. It controls
the intersection length of a translation axis of an element of a
discrete group of automorphisms of a tree with its images under this
group. We will use it in Section \ref{sec:hall} in order to prove Theorem
\ref{theo:hallfunctionfield}.

Recall that an automorphism
  $\ga$ of a simplicial tree $\TT$ with geometric realisation $|\TT|$
  is {\it loxodromic} if it fixes no point of $|\TT|$, that its {\it
    translation length} $\ell(\ga)=\min_{x\in V\TT} d(x,\ga x)$ is
  then positive and that its {\it translation axis}
  $\Ax_\ga=\{x\in|\TT|\;:\; d(x,\ga x)=\ell(\ga)\}$ is then a geodesic
  line in $|\TT|$.

\blemm \label{lem:fellowparallel} Let $\Ga$ be a discrete group of
automorphisms of a locally finite tree $\TT$. Let $\ga_0\in\Ga$ be a
loxodromic element on $\TT$. 
Let $k_0=\min_{x\in \Ax_{\ga_0}}|\Ga_x|$ be the
minimal order of the stabiliser in $\Ga$ of a vertex of $\Ax_{\ga_0}$
and let $\Ga_0$ be the stabiliser of $\Ax_{\ga_0}$ in $\Ga$. Then for
every $\ga\in \Ga-\Ga_0$, the length of the geodesic segment
$\ga\Ax_{\ga_0}\cap \Ax_{\ga_0}$ is less than $(k_0+1)\ell(\ga_0)-1$.
\elemm

\dem 
Assume for a contradiction that the length $L\in\NN$ of $\ga\Ax_{\ga_0}
\cap \Ax_{\ga_0}$ is at least $(k_0+1)\ell(\ga_0)-1$.  Denote by
$\cli xy$ the geodesic segment $\ga\Ax_{\ga_0}\cap \Ax_{\ga_0}$, such
that $\ga_0 x$ and $y$ are on same side of $x$ on $\Ax_{\ga_0}$.  Let
$\epsilon=1$ if $\ga\ga_0\ga^{-1} x$ and $y$ are on same side of $x$
on $\ga\Ax_{\ga_0}$, and $\epsilon=-1$ otherwise.

Since $\ga_0$ acts by a translation of length $\ell(\ga_0)$ on
$\Ax_{\ga_0}$, there exists a point $x'\in\cli xy$ at distance at most
$\ell(\ga_0)-1$ from $x$ such that $|\Ga_{x'}|=k_0$.  Note that $\ga
\ga_0^{\;\epsilon\,}\ga^{-1}$ acts by a translation of length
$\ell(\ga_0)$ on $\ga\Ax_{\ga_0}$ and that the translation directions
of $\ga_0$ and $\ga\ga_0^{\;\epsilon\,}\ga^{-1}$ coincide on $\cli xy$.
Hence for every $k \in \{0,1,\dots,k_0\}$, the point $\ga
\ga_0^{\;\epsilon k\,} \ga^{-1}x'$ belongs to $\cli xy$ by the
assumption on $L$, and $\ga_0^{-k} \ga \ga_0^{\;\epsilon k\,} \ga^{-1}x' 
= x'$.  Since the stabiliser of $x'$ has order less than $k_0+1$,
there are hence distinct $k,k'\in \{0,1,\dots,k_0\}$ such that
$$
\ga_0^{-k} \ga \ga_0^{\;\epsilon k\,}\ga^{-1}= \ga_0^{-k'} \ga
\ga_0^{\;\epsilon k'\,}\ga^{-1}\;,
$$
that is, $\ga_0^{k'-k} \ga = \ga\ga_0^{\epsilon (k'-k)}$.  Since
$\Ax_{\ga_0^m}=\Ax_{\ga_0}$ and $\ga\Ax_{\ga_0^{\;m'}}=
\Ax_{\ga\ga_0^{\;m'}\ga^{-1}}$ for all $m,m'\in\ZZ-\{0\}$, this
implies that $\ga$ preserves $\Ax_{\ga_0}$, a contradiction since
$\ga\notin\Ga_0$.  
\cqfd

\section{Quadratic Diophantine approximation in completions of 
function fields}
\label{sec:hall}

Let $K$ be a function field over $\FF_q$, let $v$ be a (normalised
discrete) valuation of $K$, let $R_v$ be the affine function ring
associated with $v$, and let $\Ga$ be a finite index subgroup of
$\Ga_v=\PGL_2(R_v)$ (for instance a congruence subgroup).

We denote by
$$
\QI_v = \{x \in K_v\;: \;[K(x) : K] = 2\}
$$ 
the set of quadratic irrationals in $K_v$ over $K$, and we fix
$\alpha \in\QI_v$. We denote by $\alpha^\sigma\in\QI_v$ the Galois
conjugate of $\alpha$ over $K$, and by
$$
\Theta_{\alpha,\,\Ga} =\Ga\cdot\{\alpha,\alpha^\sigma\}
$$ the union of the orbits of $\alpha$ and $\alpha^\sigma$ under the
projective action of $\Ga$, with $\Theta_{\alpha} =
\Theta_{\alpha,\,\Ga_v}$.  Note that $\alpha^\sigma\neq\alpha$, since
an irreducible quadratic polynomial over $K$ which is inseparable does
not split over $K_v$ (see for instance \cite[Lem.~17.2]{BroParPau18}),
and that there exists a loxodromic element $\ga_\alpha\in\Ga_v$ such
that $\oi\alpha{\alpha^\sigma}\;=\Ax_{\ga_\alpha}$ (see for instance
\cite[Prop.~17.3]{BroParPau18}).  For all $x \in K_v$ and
$\beta\in\QI_v$ with $x\neq \beta$, let
$$
c(x,\beta)=\frac{|x-\beta|_v}{|\beta-\beta^\sigma|_v}\in q_v^\ZZ\;.
$$

The following result gives a geometric interpretation to this quantity. 

\blemm\label{lem:calculcxbeta}
Let $x\in K_v$ and $\beta\in \QI_v$ with $x\neq \beta,\beta^\sigma$.

\smallskip\noindent 
(1) If the geodesic lines $\oi\infty x$ and
$\oi\beta{\beta^\sigma}$ in $\TT_v$ are disjoint or meet at only one
vertex, then, with $n$ the distance between them,
$$
c(x,\beta) = c(x,\beta^\sigma)= q_v^n\geq 1\;.
$$

\smallskip\noindent 
(2) If the geodesic lines $\oi\infty x$ and $\oi\beta{\beta^\sigma}$
in $\TT_v$ meet along a geodesic segment of length $n>0$, with the
closest point to $\beta$ on $\oi\infty x$ closer to $x$ than the
closest point to $\beta^\sigma$ on $\oi\infty x$, then
$$
c(x,\beta)= q_v^{-n}<1=c(x,\beta^\sigma)\;.
$$ 
In particular, $\min\{c(x,\beta), c(x,\beta^\sigma)\} < 1$.
\begin{center}
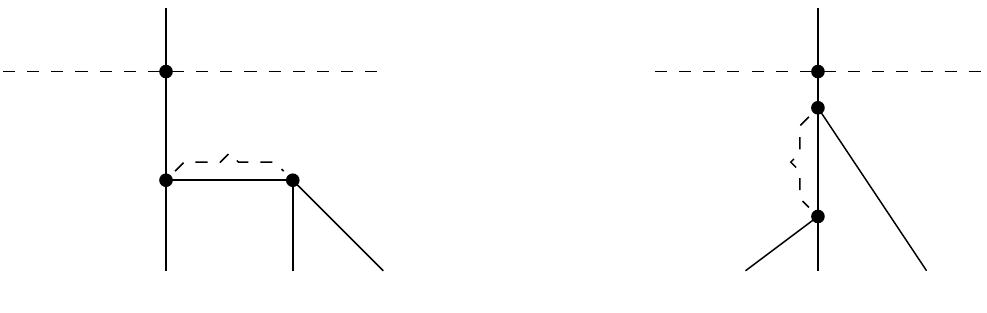
\end{center}
\elemm

\medskip
\dem 
For all distinct $y,z\in K_v$ and for every horoball $\H$ in
$\TT_v$ centered at $\infty$, let
$$
d_\H(y,z)=\lim_{t\ra+\infty} e^{\frac{1}{2}d(y_t,\,z_t)-t}\;,
$$ 
where $t\mapsto y_t$, $t\mapsto z_t$ are the geodesic lines starting
from $\infty$, through $\partial\H$ at time $t=0$, ending at the
points at infinity $y,z$ respectively.  By for instance
\cite[Eq.~(15.2)]{BroParPau18}, we have
$$
|y-z|_v=d_{\H_\infty}(y,z)^{\ln q_v}\;.
$$
Since the ratio $\frac{d_\H(y,\,z)}{d_\H(y',\,z')}$ does not depend on the
horoball $\H$ centered at $\infty$, if $\H$ is small enough, we hence
have
$$
c(x,\beta)=\frac{|x-\beta|_v}{|\beta-\beta^\sigma|_v}=
q_v^{-d(\H,\;\oi x{\,\beta})+d(\H,\;\oi\beta{\,\beta^\sigma})}\;.
$$
The result follows.
\cqfd

\medskip
For every $x \in K_v - (K\cup \Theta_{\alpha,\,\Ga})$,
we define the {\it approximation constant} of $x$ by the (extended)
$\Ga$-orbit of $\alpha$ as
$$
c_{\alpha,\,\Ga}(x)=\liminf_{\beta\in\Theta_{\alpha,\,\Ga},\;|\beta-\beta^\sigma|_v\ra 0}
\;c(x,\beta)\;.
$$ 
When $x$ is itself a quadratic irrational, the following result
gives a geometric computation of the approximation constant
$c_{\alpha,\,\Ga}(x)$.

\brema\label{rem:computquadrat}
For all $\alpha,\beta\in\QI_v$ such that $\beta\notin
\Theta_{\alpha,\,\Ga}$, we have
$$
c_{\beta,\,\Ga}(\alpha)=c_{\alpha,\,\Ga}(\beta)=q^{-n(\alpha,\,\beta)}
$$ 
where $n(\alpha,\beta)=\max_{\ga\in\Ga} \length\big(\,
\oi{\alpha^\sigma}\alpha\cap \ga \;\oi{\beta^\sigma}\beta\,\big)$ if
some image of $\oi{\beta^\sigma}\beta$ by an element of $\Ga$ meets
$\oi{\alpha^\sigma}\alpha$, and $n(\alpha,\beta)=- \min_{\ga\in\Ga}
d\big(\,\oi{\alpha^\sigma}\alpha, \ga \;\oi{\beta^\sigma}\beta\,\big)$
otherwise.  
\erema

\dem Note that since the elements of $\Ga$ preserve the lengths of
segments, and by a change of variable $\ga\mapsto\ga^{-1}$, we have
$n(\beta,\alpha)=n(\alpha,\beta)$, so that the first equality follows
from the second one.

By a proof similar to the one of Proposition \ref{lem:fellowparallel},
for all $\alpha,\beta\in\QI_v$ such that $\beta\notin
\Theta_{\alpha,\,\Ga}$, there exists a constant $\ell_{\alpha,\beta}$
(depending only on the translation lengths of primitive elements in
$\Ga$ preserving the geodesic lines $\oi{\alpha^\sigma}\alpha$ and
$\oi{\beta^\sigma}\beta$, as well as on the cardinalities of the
stabilisers in $\Ga$ of the vertices of these geodesic lines) such
that for every $\ga\in\Ga$, the length of the intersection
$\oi{\alpha^\sigma}\alpha\cap \ga \;\oi{\beta^\sigma}\beta$ is at most
$\ell_{\alpha,\beta}$.

First assume that some image of $\oi{\beta^\sigma}\beta$ by an element
of $\Ga$ meets $\oi{\alpha^\sigma}\alpha$. Using Lemma
\ref{lem:calculcxbeta} for the second equality, the fact that
$\oi{\infty}{\,\beta}$ and $\oi{\beta^\sigma}{\,\beta}$ share a subray
with endpoint $\beta$ for the third equality, and Lemma
\ref{lem:calculcxbeta} (2) for the fourth one, we have
\begin{align*}
c_{\alpha,\,\Ga}(\beta) & =
\liminf_{\alpha'\in\Theta_{\alpha,\,\Ga},\;|{\alpha'}^\sigma-\alpha'|_v\ra 0} 
c(\beta,\alpha')
=\liminf_{\alpha'\in\Theta_{\alpha,\,\Ga},\;|{\alpha'}^\sigma-\alpha'|_v\ra 0,
\;\oi{{\alpha'}^\sigma}{\,\alpha'}\;\cap\;\oi{\infty}{\,\beta}\neq 
\emptyset} c(\beta,\alpha')
\end{align*}
\begin{align*}
\\ & = \liminf_{\alpha'\in\Theta_{\alpha,\,\Ga},\;|{\alpha'}^\sigma-\alpha'|_v\ra 0,
\;\oi{{\alpha'}^\sigma}{\,\alpha'}\;\cap\;\oi{\beta^\sigma}{\,\beta}\neq 
\emptyset} c(\beta,\alpha') \\ & = 
\liminf_{{\alpha'}\in\Theta_{\alpha,\,\Ga},\;|{\alpha'}^\sigma-{\alpha'}|_v\ra 0,
\;\oi{{\alpha'}^\sigma}{\,\alpha'}\;\cap\;\oi{\beta^\sigma}{\,\beta}\neq 
\emptyset} q^{-\length(\oi{{\alpha'}^\sigma}{\,\alpha'}\;\cap
 \;\oi{\beta^\sigma}{\,\beta})}
= q^{-n(\beta,\,\alpha)}\;. 
\end{align*}

Otherwise, the result follows by using Lemma \ref{lem:calculcxbeta}
(1).  \cqfd

\medskip
We define the {\it quadratic Lagrange spectrum} of $\alpha$
relative to $\Ga$ as
$$
\Sp(\alpha,\,\Ga) = 
\{c_{\alpha,\,\Ga}(x)\;:\; x\in K_v - (K\cup \Theta_{\alpha,\,\Ga})\}\;,
$$ 
and $\Sp(\alpha) =\Sp(\alpha,\,\Ga_v)$. Note that
$\Sp(\alpha,\,\Ga)$ is contained in $q_v^{\,\ZZ}\cup\{0,+\infty\}$ and
that for every $\beta\in\Theta_{\alpha,\,\Ga}$, the functions
$c_{\alpha,\,\Ga}$ and $c_{\beta,\,\Ga}$ on $K_v - (K\cup
\Theta_{\alpha,\,\Ga})=K_v - (K\cup
\Theta_{\beta,\,\Ga})$ coincide, so that $\Sp(\alpha,\,\Ga) =
\Sp(\beta,\,\Ga)$. 

\bprop \label{prop:Lagclosborn}
The quadratic Lagrange spectrum $\Sp(\alpha,\,\Ga)$ is closed and
bounded in $\RR$.  
\eprop

The above result allows us to define the 
{\it Hurwitz constant} of $\alpha$ relative to $\Ga$ as
$$
\max \Sp(\alpha,\,\Ga)\in q_v^\ZZ\;,
$$ 
and the {\it Hurwitz constant} of $\alpha$ as $\max
\Sp(\Ga_v,\,\alpha)$.  The proof of Proposition \ref{prop:Lagclosborn}
actually gives an upper bound on $\Sp(\alpha,\,\Ga)$ which does not
depend on the quadratic irrational $\alpha$, see Equation
\eqref{eq:majounif} below. In the special case when
$(K,v,\Ga)=(\FF_q(Y), v_\infty, \Ga_{v_\infty})$, we will prove more
precisely in Section \ref{sec:comput} that $\max_{\alpha\in\QI_v}
\max\Sp(\alpha)= \frac{1}{q^2}$.

\medskip
\dem It follows from \cite[Thm.~1.6]{HerPau10}\footnote{Actually,
  Thm.~1.6 of \cite{HerPau10} is stated only for $K= \FF_q(Y)$,
  $v=v_\infty$ and $\Ga=\Ga_v$, but it has an analogous version for
  general $(K,v,\Ga)$ by using \cite[Prop.~1.5]{HerPau10}.} that if
$m_{K_v}$ is a Haar measure on the locally compact additive group of
$K_v$, then $c_{\alpha,\,\Ga}(x)=0$ for $m_{K_v}$-almost every $x\in
K_v$. Therefore $0\in \Sp(\alpha,\,\Ga)$, and the quadratic Lagrange
spectrum of $\alpha$ relative to $\Ga$ is closed.

\medskip 
Let us fix $x \in K_v - K$ and let us prove that
$c_{\alpha,\,\Ga}(x)\leq q_v^{\diam \,E_\Ga}$, where $E_\Ga$ is as
defined in Section \ref{sec:background}. This proves Proposition
\ref{prop:Lagclosborn} with a uniform bound on the Hurwitz constants
\begin{equation}\label{eq:majounif}
\forall\;\alpha\in\QI_v,\;\;\;\;\max\Sp(\alpha,\,\Ga) \leq q_v^{\diam \,E_\Ga}\;.
\end{equation}

Since $x$ is irrational and since any geodesic ray entering into a
horoball and not converging towards its point at infinity has to exit
the horoball, the geodesic line $\oi\infty x$ from $\infty$ to $x$
cannot stay after a given time in a given horoball of the family
$(\H_\xi)_{\xi\in\PP_1(K)}$ defined in Section
\ref{sec:background}.  Hence there exists a sequence $(p_n)_{n\in\NN}$
of points of $\wt E_\Ga$ converging to $x$ along the geodesic line
$\oi\infty x\,$. Since $E_\Ga=\Ga\bs\wt E_\Ga$ is finite and since no
geodesic line is contained in a horoball, there exists a sequence
$(\ga'_n)_{n\in\NN}$ in $\Ga$ such that $d(p_n,\ga'_n\,
\oi\alpha{\alpha^\sigma}) \leq \diam\; E_\Ga$ for all $n\in\NN$.

By Lemma \ref{lem:calculcxbeta}, there exists $\beta_n\in
\{\ga'_n\alpha, \ga'_n\alpha^\sigma\}\subset \Theta_{\alpha,\,\Ga}$
such that $c(x,\beta_n) < 1$ if $\oi{\beta_n}{\beta_n^\sigma}\;= \;
\ga'_n \;\oi\alpha{\alpha^\sigma}$ meets $\oi\infty x\,$ in at least
an edge, and $c(x,\beta_n) \leq q_v^{\diam \,E_\Ga}$ otherwise. Hence
$$
\liminf_{n \ra+\infty} \;c(x,\beta_n) \leq q_v^{\diam \,E_\Ga}\;.
$$ 
Let $\ga_\alpha\in\Ga_v$ be a loxodromic element  such that
$\oi\alpha{\alpha^\sigma}\;=\Ax_{\ga_\alpha}$. Since $\Ga$ has finite
index in $\Ga_v$, up to replacing $\ga_\alpha$ by a positive power, we
may assume that $\ga_\alpha$ belongs to $\Ga$. Since the length of the
intersection of two distinct translates of $\Ax_{\ga_\alpha}$ by
elements of $\Ga$ is uniformly bounded by Lemma
\ref{lem:fellowparallel}, we have $\lim_{n \ra+\infty}
\;|\beta_n-\beta_n^\sigma|_v = 0$. Hence by the definition of the
approximation constants, we have as wanted $c_{\alpha,\,\Ga}(x) \leq
q_v^{\diam \,E_\Ga}$.  \cqfd

\medskip
The following result, which implies Theorem
\ref{theo:hurwitzhallintro} (2) in the introduction, says that the
nonarchimedean quadratic Lagrange spectra contain Hall rays. Note that
its proof gives an explicit upper bound on the constant whose
existence is claimed.

\btheo\label{theo:hallfunctionfield} 
There exists $m_\alpha\in\NN$ such that for all $m\in\NN$ with $m\geq
m_\alpha$, we have $q_v^{-m}\in\Sp(\alpha,\,\Ga)$.  
\etheo

\dem Let 
$$
k_\alpha=\min_{\,x\,\in \;\oi{\alpha^\sigma}\alpha\;\cap \,V\TT_v}\; |\Ga_x|
$$ 
be the minimal order of the stabiliser in $\Ga$ of a vertex of the
geodesic line $\oi{\alpha^\sigma}\alpha$. Let $\ga_\alpha\in\Ga$ be a
loxodromic element (with minimal translation length) such that
$\Ax_{\ga_\alpha}= \oi{\alpha^\sigma}\alpha$ and let $\kappa_\alpha =
(k_\alpha+1)\ell(\ga_\alpha)-2$. By Lemma \ref{lem:fellowparallel},
for all $\beta,\beta'\in \Theta_{\alpha,\,\Ga}$, if
$\beta'\notin\{\beta,\beta^\sigma\}$, then the intersection
$\oi{\beta^\sigma}\beta\cap\oi{{\beta'}^\sigma}{\beta'}$ is a
(possibly empty) segment of length at most $\kappa_\alpha$. Take
$$
m_\alpha = 2\kappa_\alpha+1\;.
$$ 
Let us fix $m\in\NN$ with $m\geq m_\alpha$, and let us prove that
$q_v^{-m}$ belongs to $\Sp(\alpha,\,\Ga)$, which gives Theorem
\ref{theo:hallfunctionfield}. For this, let us construct $\xi\in
K_v-(K\cup\Theta_{\alpha,\,\Ga})$ and a sequence $(\beta_n)_{n\in\NN}$
in $\Theta_{\alpha,\,\Ga}$ such that 

$\bullet$~ $|\beta_n-\beta_n^\sigma|_v\ra 0$
as $n\ra+\infty$, 

$\bullet$~ the length of the intersection
$\oi{\beta_n^\sigma}{\beta_n}\, \cap\,\oi\infty\xi$ is exactly $m$,

$\bullet$~ 
the closest point to $\beta_n$ on $\oi\infty\xi$ is closer to $\xi$ than
the closest point to $\beta_n^\sigma$ on $\oi\infty\xi\,$, 

$\bullet$~ for every $\beta \in\Theta_{\alpha,\,\Ga}$, either
$\oi{\beta^\sigma}\beta$ and $\oi\infty\xi$ are disjoint or their
nonempty intersection has length at most $m$.  

\noindent By Lemma \ref{lem:calculcxbeta}, this proves that
$$
c_{\alpha,\,\Ga}(\xi)=\liminf_{\beta\in\Theta_{\alpha,\,\Ga},\;|\beta-\beta^\sigma|_v\ra 0}
\;c(\xi,\beta)=\liminf_{n\ra +\infty}
\;c(\xi,\beta_n)=q_v^{-m}\;,
$$
so that $q_v^{-m}$ does belong to $\Sp(\alpha,\,\Ga)$.

Since the image of $\oi{\alpha^\sigma}\alpha$ in $\Ga\bs\TT_v$ is
compact, and since $\Ga\cdot\infty$ is a cusp, there exists
$\beta_0\in \Theta_{\alpha,\,\Ga}$ such that if $x_0$ is the closest
point to $\infty$ on $\oi{\beta_0^\sigma}{\beta_0}$, then the open
geodesic ray $\oi{x_0}\infty$ does not meet any $\oi{\beta^\sigma}\beta$
for $\beta\in \Theta_{\alpha,\,\Ga}$. Let $y_0$ be the vertex on
$\cloi{x_0}{\beta_0}$ at distance exactly $m$ from $x_0$, and let $e_0$ be
an edge with origin $y_0$ and not contained in
$\oi{\beta_0^\sigma}{\beta_0}$.

\begin{center}
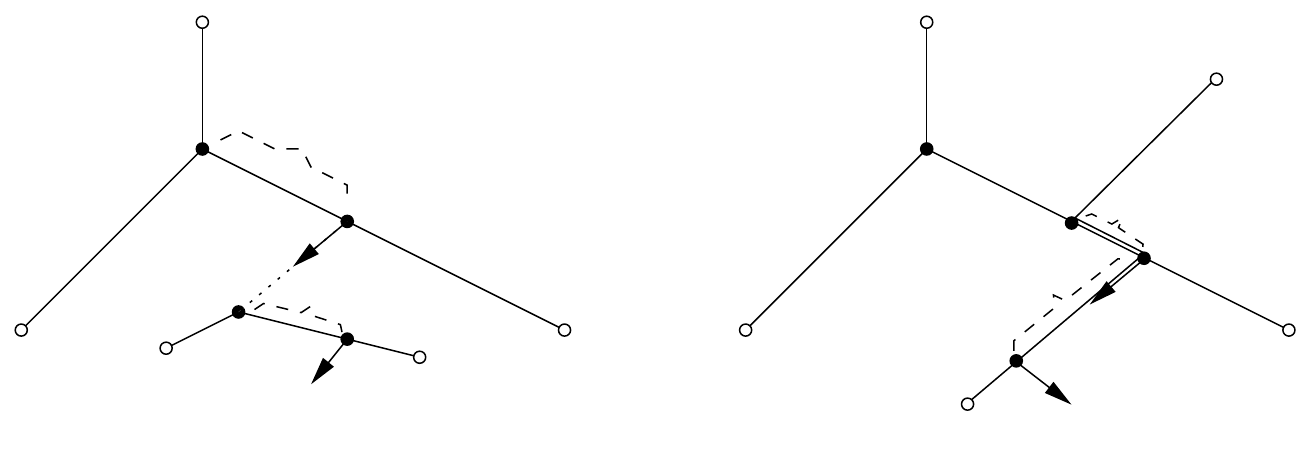
\end{center}

Assume first (see the above picture on the left) that there exists no
$\beta\in \Theta_{\alpha,\,\Ga}$ such that $e_0$ is contained in
$\oi{\beta^\sigma}\beta$. Since $\Ga$ is a lattice, the set
$\Theta_{\alpha,\,\Ga}$, which contains the orbit $\Ga\cdot\alpha$, is
dense in $\partial_\infty\TT_v$. Hence there exists $\beta_1\in
\Theta_{\alpha,\,\Ga}$ such that if $x_1$ is the closest point to
$y_0$ on $\oi{\beta^\sigma_1}{\beta_1}$, then the open segment
$\oi{y_0}{x_1}$ contains $e_0$ in its closure and meets no
$\oi{\beta^\sigma}\beta$ for $\beta\in \Theta_{\alpha,\,\Ga}$. Let
$y_1$ be the vertex on the geodesic ray $\cloi{x_1}{\beta_1}$ at
distance equal to $m$ from $x_1$, and let $e_1$ be an edge with origin
$y_1$ and not contained in $\oi{\beta_1^\sigma}{\beta_1}$. Note that
$d(y_0,y_1)\geq d(x_1,y_1) =m\geq m_\alpha\geq \kappa_\alpha+1$. By
Lemma \ref{lem:fellowparallel}, for every $\beta\in
\Theta_{\alpha,\,\Ga}$, the length of the (possibly empty)
intersection $\cloi{y_1}\infty\cap\oi {\beta^\sigma} \beta$ is at most
$m$ (exactly $m$ if $\beta\in \{\beta_0, \beta_0^\sigma,\beta_1,
\beta_1^\sigma\}$ and at most $\kappa_\alpha$ otherwise).

\medskip 
Assume now (see the above picture on the right) that there exists
$\beta_1\in \Theta_{\alpha,\,\Ga}$ such that $e_0$ is contained in
$\oi{\beta^\sigma_1}{\beta_1}$. Up to exchanging $\beta_1$ and
$\beta_1^\sigma$, we may assume that $e_0$ is contained in
$\cloi{y_0}{\beta_1}$.

Let $x_1\in V\TT_v$ be such that $\cloi{y_0}\infty\cap\cloi{y_0}
{\beta_1^\sigma}=\cli{y_0}{x_1}$. Note that by Lemma
\ref{lem:fellowparallel}, we have $d(x_1,y_0)\leq \kappa_\alpha$ and
$x_1\in \cloi{y_0}{x_0}$ since $d(x_0,y_0)=m\geq
m_\alpha>\kappa_\alpha$. Let $y_1$ be the point at distance equal to
$m-d(x_1,y_0)$ from $y_0$ on $\cloi{y_0}{\beta_1}$. Note that
$d(y_0,y_1)\geq \kappa_\alpha+1$ since $m\geq
m_\alpha=2\kappa_\alpha+1$, and in particular by Lemma
\ref{lem:fellowparallel}, there is no $\beta\in \Theta_{\alpha,\,\Ga}$
such that $\oi{\beta^\sigma}\beta$ contains both $e_1$ and $y_0$. By
Lemma \ref{lem:fellowparallel}, for every $\beta\in
\Theta_{\alpha,\,\Ga}$, the length of
$\cloi{y_1}\infty\cap\oi{\beta^\sigma}\beta$ is at most $m$ (exactly
$m$ if $\beta\in\{\beta_0, \beta_0^\sigma,\beta_1, \beta_1^\sigma\}$
and at most $2\kappa_\alpha$ otherwise).

\medskip By iterating this construction, we obtain sequences
$(\beta_n)_{n\in\NN}$ in $\Theta_{\alpha,\,\Ga}$, $(x_n)_{n\in\NN}$,
$(y_n)_{n\in\NN}$ in $V\TT_v$ and $(e_n)_{n\in\NN}$ in $E\TT_v$ such
that, for all $n \in\NN$,
\begin{itemize}
\item $y_n\in\cloi{y_{n+1}}\infty$ and $d(y_n,y_{n+1})\geq 
\kappa_\alpha+1$,

\item for every $\beta\in \Theta_{\alpha,\,\Ga}$, the length of
$\cloi{y_n}\infty\cap\oi{\beta^\sigma}\beta$ is at most $m$, and the
length of $\cloi{y_n}\infty\cap\oi{\beta^\sigma_n}{\beta_n}= \cli{y_n}{x_n}$
is exactly $m$,

\item $y_n\in\cloi{x_n}{\beta_n}$,

\item $e_n$ exits $\oi{\beta^\sigma_n}{\beta_n}$ at $y_n$ and is
contained in $\cli{y_n}{y_{n+1}}$,

\item $x_{n+1}$ and $y_n$ belong to $\cli{x_n}{y_{n+1}}$, and either
  $y_n$ belongs to $\cloi{x_n}{x_{n+1}}$ in which case
  $\oi{y_n}{x_{n+1}}$ meets no $\oi{\beta^\sigma}\beta$ for $\beta\in
  \Theta_{\alpha,\,\Ga}$, or $y_n$ belongs to $\cloi{x_{n+1}}{y_{n+1}}$.
\end{itemize}
By the first point, the sequence $(y_n)_{n\in\NN}$ converges to
$\xi\in\partial_\infty \TT_v$, such that $y_0,y_1,y_2,\dots$ are in
this order on the geodesic line $\oi\infty\xi$ oriented from $\infty$
to $\xi$. The point at infinity $\xi$ belongs neither to $K$ (since
$y_n$ belongs to $\oi{\beta^\sigma_n}{\beta_n}$ and the geodesic lines
$\oi{\beta^\sigma}\beta$ for $\beta\in\Theta_{\alpha,\,\Ga}$ do not
enter a small enough horoball centered at any point of $K$) nor to
$\Theta_{\alpha,\,\Ga}$ (otherwise $\oi\infty\xi\cap
\oi{\xi^\sigma}\xi$ would be infinitely long, contradicting the second
point for $n$ large enough).

The length of $\oi\infty\xi\cap\oi{\beta^\sigma_n}{\beta_n}$, which
is equal to $\cli{y_n}{x_n}$ since $e_n$ exits $\oi{\beta^\sigma_n}{
\beta_n}$ at $y_n$, is exactly $m$, and in particular is bounded in
$n$.  Since $(y_n)_{n\in\NN}$ converges to the point at infinity
$\xi$, we have $|\beta_n-\beta_n^\sigma|_v\ra 0$ as $n\ra+ \infty$. By
construction, there is no $\beta\in \Theta_{\alpha,\,\Ga}$ such that
the length of $\oi\infty\xi\cap\oi{\beta^\sigma}\beta$ is larger
than $m$.

Therefore $\xi$ satisfies the properties required at the beginning of
the proof, and Theorem \ref{theo:hallfunctionfield} follows.  
\cqfd

\medskip
In the next section, we will give several computations, using the
continued fraction expansions, in the special case when $K=\FF_q(Y)$,
$v=v_\infty$ is the valuation at infinity, and $\Ga=\Ga_{v_\infty}$ is
the full Nagao lattice $\PGL_2(\FF_q[Y])$.

\section{Computations of approximation constants, Hurwitz constants 
and quadratic Lagrange spectra for fields of formal Laurent series}
\label{sec:comput}

In this section, we use the notation $\FF_q$, $R=\FF_q[Y]$,
$K=\FF_q(Y)$, $\wh K=\FF_q((Y^{-1}))$, $|\cdot|$ given in the
introduction.  Let
$$
\OOO=\{f\in\wh K: |f|\leq 1\}=\FF_q[[Y^{-1}]]
$$ 
be the ring of formal power series in one variable $Y^{-1}$ over
$\FF_q$.  Its unique maximal ideal is $\mmm=Y^{-1}\OOO$. We denote by
$\TT$ the Bruhat-Tits tree of $(\PGL_2,\wh K)$, with standard base
point $*=[\OOO\times\OOO]$, and $\Ga=\PGL_2(R)$, see Section
\ref{sec:background}.

Any element $f\in \wh K$ may be uniquely written as a sum
$f=[f]+\{f\}$ of its {\em integral part} $[f]\in R$ and its {\em
  fractional part} $\{f\}\in \mmm$.  The {\it Artin map} $\Psi:
\mmm-\{0\}\ra \mmm$ is defined by $f\mapsto \big\{\frac{1}{f}\big\}$.
Any $f\in \wh K-K$ has a unique continued fraction expansion
$$
f=[a_0,a_1,a_2,a_3,\dots]=a_0 +
\cfrac{1}{a_1+\cfrac{1}{a_2+ \cfrac{1}{a_3+\ddots}}}
$$ 
with $a_0=a_0(f)=[f]\in R$, and $a_n=a_n(f)=
\big[\frac{1}{\Psi^{n-1}(f-a_0)} \big]\in R$ a nonconstant polynomial
for $n\geq 1$. The polynomials $a_n(f)$ are called the {\it
  coefficients} of the continued fraction expansion of $f$.  For every
$n\in\NN$, the rational element
$$
\frac{P_n}{Q_n}=[a_0,a_1,a_2,\dots,a_{n-a},a_n]=a_0 +
\cfrac{1}{a_1+\cfrac{1}{a_2+ \cfrac 1{{\phantom{\big|}}^{\ddots}+
\cfrac{{\phantom{\big|}}^{\ddots}\, 1}{a_{n-1}+
\cfrac{1}{a_{n}}}}}}\;
$$ 
is the $n$-th {\em convergent} of $f$. We refer to
\cite{Lasjaunias00, Schmidt00,Paulin02} for details and further
information on continued fraction expansions of formal Laurent series
and their geometric interpretation in terms of the Bruhat-Tits tree
$\TT$.

For every $f\in \wh K-K$, let 
\begin{align*}
M(f) & =\limsup_{k\to+\infty}\;\deg a_k\;\geq 1\,,\\
M_2(f) & =\limsup_{k\to+\infty}\big(\deg a_k+\deg a_{k+1}\big)\;\geq 2\,,\\
m(f) & =\liminf_{k\to+\infty}\;\deg a_k\;\geq 1\,.
\end{align*}

As recalled in the introduction, an irrational element $\alpha\in\wh K
-K$ is quadratic over $K$ if and only if its continued fraction
expansion is eventually periodic: For every $p\in\NN$ large enough,
the sequence of coefficients $(a_{k+p}(\alpha))_{k\in\NN}$ is periodic
with period $m\in\NN-\{0\}$ and, as usual, we then write the continued
fraction expansion of $\alpha$ as $\alpha= [a_0,a_1,a_2,\dots,a_{p-1},
  \overline{a_{p},a_{p+1},\dots,a_{p+m-1}}\;]$.  We then have
$$
M(\alpha) =\max_{0\le k\le m-1}\deg a_{p+k}(\alpha)\;\;\;\textrm{ and} \;\;\;
m(\alpha) = \min_{0\le k\le m-1}\deg a_{p+k}(\alpha)\;.
$$


\medskip
Let us recall from \cite[\S 6.3, Rem.~2]{Paulin02} the penetration
properties of the geodesic lines of $\TT$ inside the canonical
equivariant family of pairwise disjoint open horoballs whose closures
cover $\TT$, see the picture below.  For every $f\in\wh K-K$, the
geodesic line $\oi \infty f$ (oriented from $\infty$ to $f$) starts in
$H_0(f)=\H_\infty$ by an initial negative subray; after $\H_\infty$,
it successively passes through an infinite sequence of open horoballs
in the family $(\H_\xi)_{\xi\in\PP^1(K)}$, denoted by
$(H_n(f))_{n\in\NN-\{0\}}$ and for every $n\in\NN$, we have
$H_{n+1}(f)=\H_{P_n/Q_n}$, where $P_n/Q_n$ is the $n$-th convergent of
$f$. For every $n\in\NN$, the maximum height of a point of
$]\infty,f[$ inside $H_{n+1}(f)$ is equal to $\deg a_{n+1}$, or,
    equivalently, the intersection of $H_{n+1}(f)$ with $]\infty,f[$
    is an open segment of length $2\deg a_{n+1}$.

\begin{center}
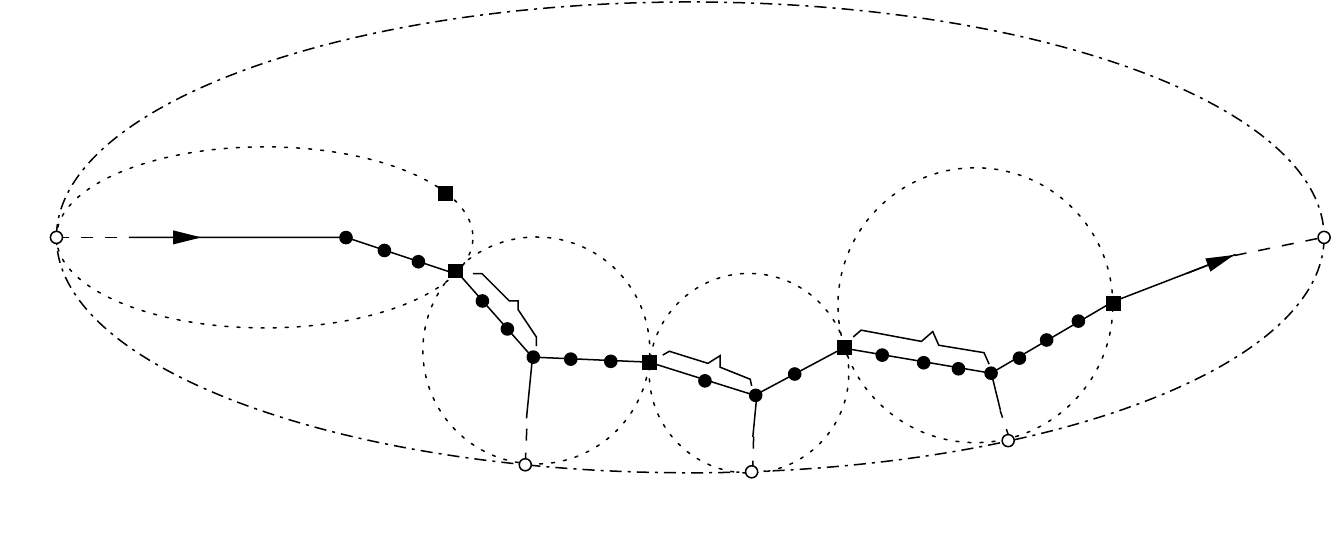
\end{center}

We will need the following lemmas, which also follow from the
geometric interpretation of the continued fraction expansion given in
\cite{Paulin02}, in order to estimate the quadratic approximation
constants of elements $f\in\wh K-K$.

\blemm\label{lem:coincidence} Let $f,f'\in\wh K-K$ with
$a_0(f)=a_0(f')$. If the geodesic lines $\oi\infty f$ and $\oi\infty
{f'}$ coincide inside the horoball $\H_{a_0(f)}$, then the polynomial
$a_1(f)-a_1(f')$ is constant.  \cqfd
\elemm

\blemm\label{lem:coincidencem}
Let $f\in\wh K-K$ and $\alpha\in \QI$ with purely periodic continued
fraction expansion. If there exist $p,q\in\NN$ and $m\in\NN$ such
that $a_{p+1}(f)= a_{q+1}(\alpha), \dots, a_{p+m}(f)=a_{q+m}(\alpha)$,
then there exists $\ga\in \Ga$ such that the geodesic lines $\oi\infty
f$ and $\ga\,\oi{\alpha^\sigma}{\alpha}$ coincide at least in $H_{p+1}(f)\cup
\dots \cup H_{p+m}(f)$, so that their intersection has length at least
$$
2\sum_{i=1}^{m}\;\deg a_{p+i}(f)\;. \;\;\;\Box
$$
\elemm

\bigskip
The proof of Proposition \ref{prop:Lagclosborn}, and in particular
Equation \eqref{eq:majounif} since $\diam \;E_\Ga =0$ for the Nagao
lattice $\Ga$ as seen in Example \ref{exem:2.2}, shows that the
quadratic Lagrange spectrum with respect to the valuation at infinity
of any quadratic irrational is contained in
$\{0\}\cup\{q^{-n}\;:\;n\in\NN\}$.  The following result, which
implies Theorem \ref{theo:hurwitzhallintro} (1) in the introduction,
improves the upper bound of the spectrum. In Corollary
\ref{coro:hurwitzconst} and Proposition
\ref{prop:hurwizconstantperiodtwo}, we will show that this upper bound
is realised for certain quadratic irrationals.

\bprop\label{prop:upperboundofspectrum} For every quadratic irrational
$\alpha$, the spectrum $\Sp(\alpha)$ is contained in
$\{0\}\cup\{q^{-n-2}\;:\;n\in\NN\}$.  
\eprop

\dem Up to replacing $\alpha$ by an element of $\Theta_\alpha$, we may
assume that $\alpha\in\mmm$ and $\alpha^\sigma\in\wh K-\OOO$, so that
the base point $*$ belongs to the geodesic line
$\oi\alpha{\alpha^\sigma}$.

Let $f\in\wh K-(K\cup \Theta_\alpha)$. Since $f$ is irrational, no
positive subray of the geodesic line $\oi\infty f$ (oriented from
$\infty$ to $f$) is contained in a single horoball of the canonical
family $(\H_\xi)_{\xi\in\PP^1(K)}$. Hence there exists a sequence
$(x_n=\ga_n *)_{n\in\NN}$, with $\ga_n\in\Ga$, of vertices in the
$\Ga$-orbit of the base point $*$, converging to $f$ on the geodesic
line $\oi\infty f$. As $*$ belongs to $\oi\alpha{\alpha^\sigma}$, we
thus have that $\oi\infty f$ and $\oi{\ga_n\cdot\alpha}{\ga_n\cdot
  \alpha^\sigma}$ meet at least at the vertex $\ga_n*$.

The stabiliser $\PGL_2(\FF_q)$ of $*$ in $\Ga$ acts transitively on
the set of pairs of distinct elements of the link of $*$. Thus, up to
multiplying $\ga_n$ on the right by an element of $\PGL_2(\FF_q)$, the
geodesic line $\oi\infty f$ meets $\oi{\ga_n\cdot\alpha} {\ga_n\cdot
  \alpha^\sigma}$ in a segment of length at least $2$ for all $n\in
\NN$.  Thus, by Lemma \ref{lem:calculcxbeta} (2), we have 
$$
\min \{c(x,\ga_n\cdot\alpha), c(x,\ga_n\cdot\alpha^\sigma)\} \leq q^{-2}\;.
$$ 
Since $|\ga_n\cdot\alpha-\ga_n\cdot\alpha^\sigma|$ tends to $0$ by
Lemma \ref{lem:fellowparallel}, we have 
$$c_\alpha(f)=
\liminf_{\beta\in\Theta_{\alpha},\;|\beta^\sigma-\beta|\ra0} c(x,\beta) \leq
q^{-2}\;.
$$ Hence the result follows. 
\cqfd

\medskip
We are now going to give a series of computations of quadratic
approximation constants. We start by two preliminary results.

\blemm\label{lem:M2} 
Let $\alpha\in\QI$ and let $f\in\wh K -(K\cup\Theta_\alpha)$.

\smallskip \noindent 
(1) If $m(f)>M(\alpha)$, then $c_\alpha(f)\ge q^{-M_2(\alpha)}$.

\smallskip\noindent 
(2) If $M(f)<m(\alpha)$, then $c_\alpha(f)\ge q^{-M_2(f)}$. 
\elemm 

\dem (1) By the penetration properties of geodesic lines in the
horoballs of the canonical family $(\H_\xi)_{\xi\in\PP^1(K)}$, for
every $\beta\in\Theta_\alpha$, the maximum height the geodesic line
$\oi{\beta^\sigma}\beta$ enters in one of these horoballs is
$M(\alpha)$. Similarly, the minimum height the geodesic line
$\oi{\infty}f$ enters one of these horoballs except finitely many of
them is $m(f)$, which is strictly bigger than $M(\alpha)$.  Hence for
all $\beta\in\Theta_\alpha$, the geodesic lines $\oi\infty f$ and
$\oi{\beta^\sigma}\beta$ can meet at most in two consecutive horoballs
$H_n(f)$ for $n\in\NN$ large enough, and their intersection has length
at most $M_2(\alpha)$ (and even at most $M(\alpha)\leq M_2(\alpha)$ if
$\oi{\beta^\sigma}\beta$ meets at most one of the horoballs $H_n(f)$
for $n\in\NN$ large enough), see the picture below.

\begin{center}
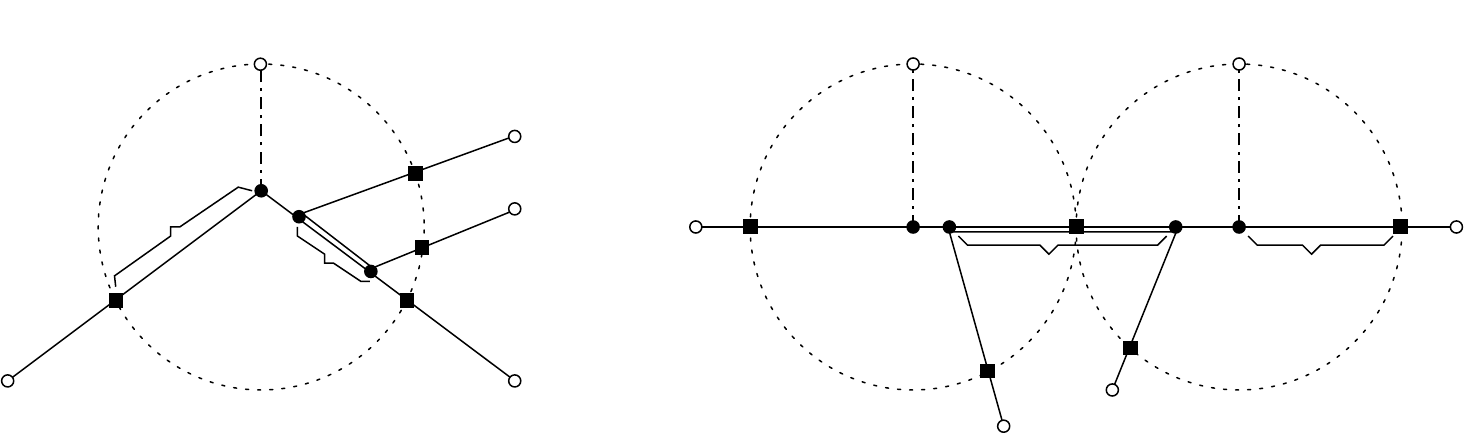
\end{center}

By Lemma \ref{lem:calculcxbeta}, we have $c(f,\beta)\ge
q^{-M_2(\alpha)}$ for all $\beta\in\Theta_\alpha$, which proves the
result.

\medskip The proof of Assertion (2) is similar.
\cqfd

\blemm\label{lem:inspectrum} For all $\alpha\in\QI$ and
$f\in\wh K-(K\cup\Theta_\alpha)$ such that $m(f)>M(\alpha)$, there
exists a sequence $(\beta_n)_{n\in\NN}$ in $\Theta_\alpha$ such that
$c(f,\beta_n)=q^{- M_2(\alpha)}$ and $|\beta_n-\beta_n^\sigma| \to 0$
as $n\to+\infty$.  
\elemm

\dem Replacing $\alpha$ by an element in its $\Ga$-orbit if necessary,
we can assume that the continued fraction expansion of $\alpha$ is
periodic, that $\alpha\in\mmm$ and $\alpha^\sigma\in\wh K-\OOO$, and
that $M_2(\alpha)=\deg a_1(\alpha)+\deg a_2(\alpha)$.  The unipotent
upper triangular subgroup $\Big\{\begin{pmatrix}1 & x\\ 0 & 1 
\end{pmatrix}\;:\;x\in R\Big\}$ of $\Ga$ fixes $\infty\in
\partial_\infty \TT$, and acts transitively on the subset $R$ of
$\partial_\infty\TT$. Since the horoballs in the canonical family
$(\H_\xi)_{\xi\in\PP^1(K)}$ whose closure meets the closure of
$\H_\infty$ are (besides $\H_\infty$ itself) the ones centred at an
element of $R$, the group $\Ga$ acts transitively on the ordered pairs
of horoballs in this family whose closures meet at one point. In
particular, for all $n\in\NN$ large enough, there exists $\ga_n\in\Ga$
sending $H_1(\alpha)$ to $H_{n}(f)$ and $H_2(\alpha)$ to $H_{n+1}(f)$.

\begin{center}
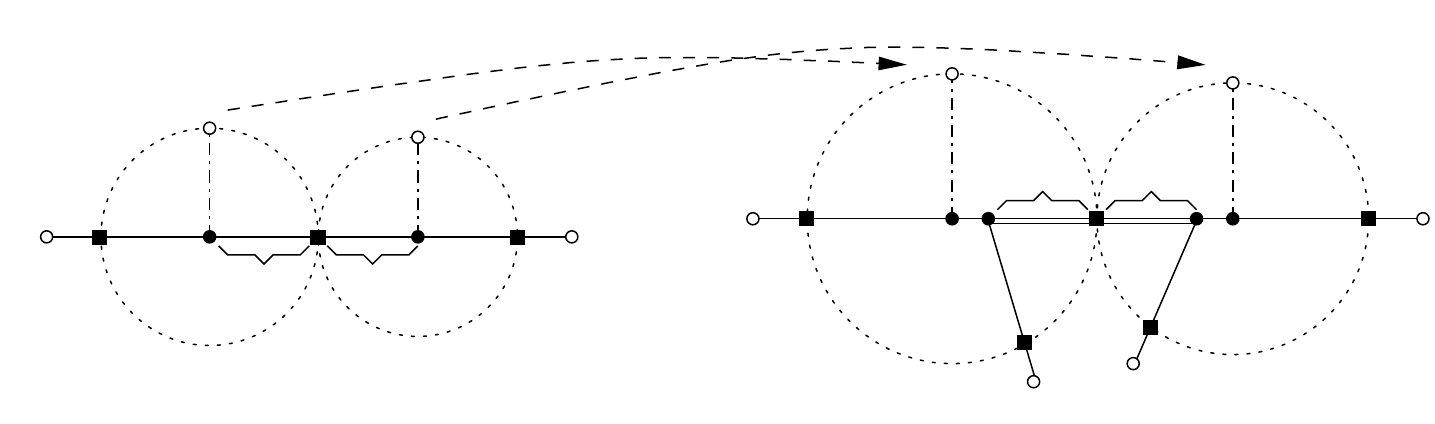
\end{center}

For all $x\in\TT$ and $\ga\in\Aut(\TT)$, and for every horoball $\H$
in $\TT$, the height of $\ga x$ with respect to $\ga\H$ is equal to
the height of $x$ with respect to $\H$. Hence for $n$ large enough,
since $m(f)>M(\alpha)$, the geodesic lines $\oi\infty f$ and
$\ga_n\,\oi{\alpha^\sigma} \alpha= \oi{\ga_n\cdot\alpha^\sigma}
{\ga_n\cdot\alpha}$ intersect along a segment of length equal to $\deg
a_1(\alpha)+\deg a_2(\alpha)=M_2(\alpha)$.  Let
$\beta_n=\ga_n\cdot\alpha$, we then have $c(f,\beta_n)=q^{-
  M_2(\alpha)}$ by Lemma \ref{lem:calculcxbeta}.  The fact that
$|\beta_n-\beta_n^\sigma| \to 0$ as $n\to+\infty$ follows from Lemma
\ref{lem:fellowparallel}. \cqfd

\bcoro\label{coro:inspectrum} Let $\alpha\in\QI$. If $f\in\wh
K-(K\cup\Theta_\alpha)$ satisfies $m(f)>M(\alpha)$, then
$c_\alpha(f)=q^{-M_2(\alpha)}\in\Sp(\alpha)$.
\ecoro
\medskip

\dem
This follows immediately from Lemmas \ref{lem:M2} (1) and
\ref{lem:inspectrum}, by the definition of the quadratic approximation
constants.
\qed

\bcoro\label{coro:hurwitzconst} Let  $\alpha\in\QI$. Then

\smallskip\noindent (1) $\max\Sp(\alpha)\ge q^{-M_2(\alpha)}$.

\smallskip\noindent (2) If $M(\alpha)=1$ or $m(\alpha)\ge 2$, then 
$\max\Sp(\alpha)=\max_{P\in\FF_q[X],\;\deg P=1} c_\alpha([\;\ov{P}\;])=q^{-2}$.  
\ecoro

\dem (1) This follows from Corollary \ref{coro:inspectrum} since for instance
$$
f=[\;\overline{a_{p+1}(\alpha)+Y^{M(\alpha)+1}, a_p(\alpha)+Y^{M(\alpha)+1}}\;]
$$ 
satisfies the assumption of Corollary \ref{coro:inspectrum} if
$p\in\NN$ is such that $M_2(\alpha)=\deg a_p(\alpha)+\deg
a_{p+1}(\alpha)$.

\medskip (2) If $M(\alpha)=1$, then $M_2(\alpha)=2$ and this follows
from Proposition \ref{prop:upperboundofspectrum} and Corollary
\ref{coro:inspectrum}. If $m(\alpha)\ge 2$, with $f=[\,\ov{Y}\,]$, we
have $M(f)=1<m(\alpha)$ and $M_2(f)=2$, thus $c_\alpha(f)\geq q^{-2}$ by Lemma
\ref{lem:M2} (2), and the result follows from Proposition
\ref{prop:upperboundofspectrum}.  
\qed

\medskip
The above corollary shows that the maximum Hurwitz constant is
attained for many quadratic irrationals $\alpha$. In fact, the same
holds for all quadratic irrationals with small enough period length.

\bprop\label{prop:hurwizconstantperiodtwo} If $\alpha$ is a quadratic
irrational over $K$ in $\wh K$ whose period of its continued fraction
expansion contains at most $q-2$ coefficients of degree $1$, then
$\max\Sp(\alpha)=\max_{P\in\FF_q[X],\;\deg P=1} c_\alpha([\;\ov{P}\;])=q^{-2}$.  
\eprop

The first equality is a (strengthened) nonarchimedean version of the
$2$-periodic case of Bugeaud's conjecture solved by Lin
\cite[Rem.~1.3]{Lin18}. In particular, if $\alpha\in\QI$ is eventually
$k$-periodic with $k\leq q-1$ (as for instance with $k=2$ since $q\geq
3$), then $\max\Sp(\alpha)=\max_{P\in\FF_q[X],\;\deg P=1}
c_\alpha([\,\ov{P}\,])=q^{-2}$. Indeed, either all coefficients of the
period of $\alpha$ have degrees $1$, in which case $M(\alpha)=1$ and
Corollary \ref{coro:hurwitzconst} (2) applies, or $\alpha$ satisfies
the assumption of Proposition \ref{prop:hurwizconstantperiodtwo}. This
proves Theorem \ref{theo:maxhurwitzintro} in the introduction.

\medskip
\dem 
Since $\card(\FF_q-\{0\})=q-1$ and by the assumption, there exists a
polynomial $P\in R$ of degree $1$ such that for every degree $1$
coefficient $a_i(\alpha)$ of the period of $\alpha$, the polynomial
$P-a_i(\alpha)$ is nonconstant. Let $f=[\,\overline P\,]$.  For all
$\beta\in\Theta_\alpha$ and $n\in\NN$ large enough, we claim that
$\oi{\beta^\sigma}\beta$ agrees with $\oi\infty f$ on a segment with
length at most $1$ inside any horoball $H_n(f)$ for $n\in\NN$. By an
argument as in the proof of Proposition
\ref{prop:upperboundofspectrum}, this implies that $c_\alpha(f)=
q^{-2}$.  This in turn implies that $\max\Sp(\alpha)\geq q^{-2}$, and
the result follows since $q^{-2}$ is an upper bound on
$\max\Sp(\alpha)$ (see Proposition \ref{prop:upperboundofspectrum}).

Assume for a contradiction that the geodesic segment
$\oi{\beta^\sigma}\beta$ agrees with $\oi\infty f$ on a segment of
length at least $2$ inside $H_n(f)$. Since $\deg a_n(f)=\deg P=1$,
this implies that $\oi{\beta^\sigma}\beta$ and $\oi\infty f$ actually
coincide inside of $H_n(f)$. Assume that the orientations of the
geodesic lines $\oi{\beta^\sigma}\beta$ and $\oi\infty f$ respectively
from $\beta^\sigma$ to $\beta$ and from $\infty$ to $f$ agree. By
Lemma \ref{lem:coincidence}, this implies that if $a_i(\beta)$ is the
coefficient in the period of $\beta$ such that $H_n(f)=H_{i}(\beta)$,
then the polynomial $P - a_i(\beta)$ is constant.  This implies that
$\deg a_i(\beta)=1$ and this contradicts the definition of $P$, since
$\alpha$ and $\beta$ have the same period (up to a cyclic
permutation).  \cqfd

\medskip
If the period of a quadratic irrational $\alpha$ is longer than $q-1$,
then its Hurwitz constant $\max\Sp(\alpha)$ may be arbitrarily small,
as the following result shows.

\bprop\label{prop:hurwitzupperbound} For all $m,k\in\NN-\{0,1\}$, let
us denote by $\{b_1,\dots, b_N\}$ the set\footnote{with an arbitrary
  ordering and $N=(q^{k+1}-q)^m$} of the finite sequences of length
$m$ of polynomials of degree at least $1$ and at most $k$, and let
$$
\E_{m,\,k}=\bigcup_{\sigma\in{\mathcal S}_N} \Ga\cdot
[\;\overline{b_{\sigma(1)}, \dots, b_{\sigma(N)}}\;]\;.
$$ 
For every $\alpha\in \E_{m,\,k}$, we have $\max\Sp(\alpha)\le
q^{-\min(2m,\,k+2)}$.  \eprop

\dem Let us fix $\alpha\in \E_{m,\,k}$ and $f\in\wh K-(K\cup
\Theta_\alpha)$, and let us prove that $c_\alpha(f)\le
q^{-\min(2m,\,k+2)}$, which gives the result.

Assume first that there exists a subsequence of coefficients
$(a_{i_n}(f))_{n\in\NN}$ of $f$ such that $i_0\geq 1$ and $\deg
a_{i_n}(f)>k$ for all $n\in \NN$. Let $x_n$ be the point of $\partial
H_{i_n}(f)$ at which the geodesic line $\oi{\infty}f$ exits the
horoball $H_{i_n}(f)$.  Since any polynomial $P\in R$ of degree $k$
occurs as a coefficient of $\alpha$ (in its periodic part), there
exists a horoball $\H_k$ in the canonical family
$(\H_{\xi})_{\xi\in\PP^1(K)}$ which intersects the geodesic
$\oi{\alpha^\sigma} \alpha$ in a segment of length exactly $2k$. Let
$y$ be the point of $\partial \H_k$ at which the geodesic
$\oi{\alpha^\sigma}\alpha$ exits the horoball $\H_k$.

\begin{center}
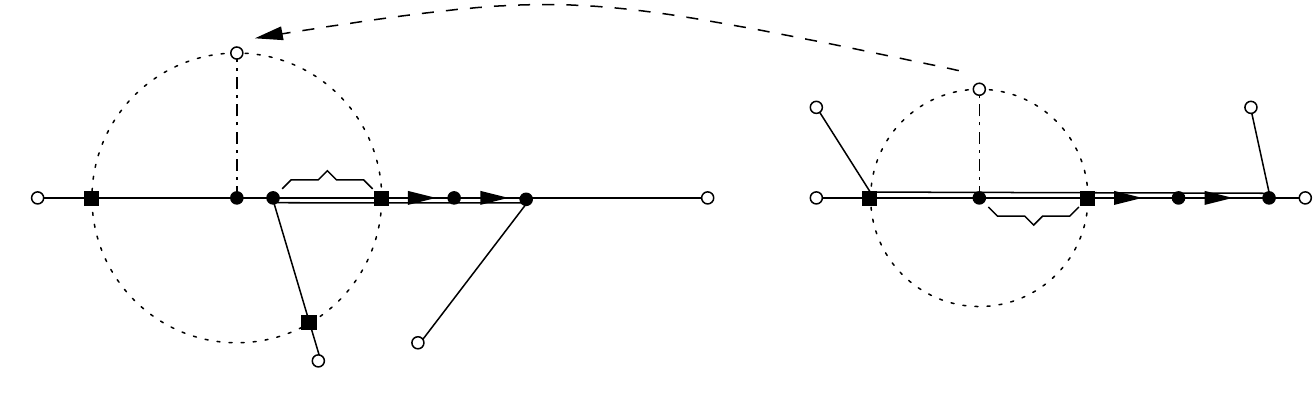
\end{center}

By the transitivity of the action of $\Ga$ on the pairs $(\H,x)$ where
$\H$ belongs to the canonical family $(\H_{\xi})_{\xi\in\PP^1(K)}$ and
$x$ belongs to $\partial \H$, there exists $\ga_n\in\Ga$ which sends
$(\H_k,y)$ to $(H_{i_n}(f),x_n)$. Let $e_1,e_2$ be the consecutive
edges along $\oi{\infty}f$ such that the origin of $e_1$ is $x_n$.
Note that the endpoint of $e_2$ (which is at distance $2$ from $x_n$)
might belong to the orbit $\Ga *$ of the base point of $\TT$ or
not. Since any pair $(P,P')$ with $P,P'\in R$, $\deg P=k$ and $\deg
P'\in\{1,2\}$ occurs as a pair of consecutive coefficients in the
continued fraction expansion of $\alpha$, there exists
$\beta_n\in\Theta_\alpha$ such that the geodesic line
$\oi{\beta_n^\sigma}{\beta_n}$ contains $\ga_n^{\,-1} e_1$ and
$\ga_n^{\,-1}e_2$, and coincides with $\oi{\alpha^\sigma} \alpha$
inside $\H_k$. Hence $\ga_n \,\oi{\beta_n^\sigma}{\beta_n}$ coincides
with $\oi{\infty}f$ on $e_1,e_2$ as well as on a segment of length
exactly $k$ inside $H_{i_n}(f)$ (since $\deg a_{i_n}(f)>k$). Thus
$$
\min \{c(f,\ga_n\cdot\beta_n^\sigma),
\;c(f,\ga_n\cdot\beta_n)\} \leq q^{-(k+2)}
$$ 
by Lemma \ref{lem:calculcxbeta}. Since
$|\ga_n\cdot\beta_n-\ga_n\cdot\beta_n^\sigma| \to 0$ as $n\to+\infty$
by Lemma \ref{lem:fellowparallel}, we have $c_\alpha(f)\le
q^{-(k+2)}$.

Otherwise, there exists $i_0\in\NN$ such that $\deg a_i(f)\le k$ for
all $i\ge i_0$.  For every $n\in\NN$, consider the string of $m$
consecutive horoballs $H_{i_0+nm}(f),\dots,H_{i_0+(n+1)m-1}(f)$
crossed by the geodesic line $\oi\infty f$. The quadratic irrational
$\alpha$ has been constructed in such a way that there exists a string
of consecutive coefficients in the period of $\alpha$ that agrees with
those of $f$ inside the above horoballs.  Using Lemma
\ref{lem:coincidencem}, this gives the estimate $c_\alpha(f)\le
q^{-2m}$. Together with the previous estimate, this completes the
proof.  
\cqfd

\bexem If $\alpha$ belongs to the set $\E_{m,k}$ constructed in
Proposition \ref{prop:hurwitzupperbound} with $k=m=2$, then
$$ \max\Sp(\alpha)=\max_{P\in R} c_\alpha([\;\ov{P}\;])=
\limsup_{P\in  R,\;\deg P\ra+\infty} c_\alpha([\;\ov{P}\;]) = q^{-4}\;.
$$ 
Indeed, Proposition \ref{prop:hurwitzupperbound} gives
$\max\Sp(\alpha)\le q^{-4}$. For any $P\in R$ with $\deg P\geq 3$ and
$f=[\,\overline P\,]$, we have $c_\alpha(f)= q^{-4}$ by Corollary
\ref{coro:inspectrum}. The result follows.  
\eexem

\medskip
Let 
$$
\varphi=[\,\overline Y\,]\,. 
$$ 
As a special case of Corollary \ref{coro:hurwitzconst}, we get
$\max\Sp(\varphi)=q^{-2}$. The following result gives in particular a
complete description of the quadratic Lagrange spectrum of $\varphi$,
proving Theorem \ref{theo:goldenratiointro} in the introduction.

\btheo\label{theo:goldenspec}
For every $P\in R$ with $\deg P=1$, we have 
$$
\Sp([\,\overline{P}\,])=\{0\}\cup\{q^{-(2+n)}:n\in\NN\}\;.
$$
\etheo

\dem Let $\alpha=[\,\overline{P}\,]$. For every $m\in\NN$, let
$f_m=[0,\overline{Y^2,P,P,\dots, P}\,]$ and $f'_m=[0,
  \overline{Y^2,P,P,\dots, P, P+1}\,]$, which are elements of
$K^{(2)}-\Theta_\alpha$ with period lengths of the periodic part of
their continued fraction equal to $m+1$ and $m+2$ respectively. Let us
prove that $c_\alpha(f_m)=q^{-2m -2}$ and $c_\alpha(f'_m)= q^{-2m
  -3}$, which implies that $q^{-n -2}$ belongs to $\Sp(\alpha)$ for
all $n\in\NN$. This gives the result by Proposition
\ref{prop:upperboundofspectrum}.

We know that $c_\alpha(f_0)=q^{-2}$ by Corollary
\ref{coro:inspectrum}.  Assume that $m\geq 1$. Note that for every
$n\in\NN$, the maximum height of the geodesic line $\oi{\infty}{f_m}$
inside the horoball $H_{n(m+1)+i+1}(f_m)$ is equal to $2$ if $i=0,m+1$ and
to $1$ if $1\leq i\leq m$. Since the geodesic lines
$\oi{f_m^\sigma}{f_m}$ and $\oi{\alpha^\sigma}{\alpha}$ both contain
points in $\Ga *$, some image of $\oi{\alpha^\sigma}{\alpha}$ by an
element of $\Ga$ meets $\oi{f_m^\sigma}{f_m}$. Hence by Remark
\ref{rem:computquadrat}, we have $c_\alpha(f_m)=q^{-n(f_m,\alpha)}$
where
$$
n(f_m,\alpha)= \max_{\ga\in\Ga}\length\big(\,
\oi{f_m^\sigma}{f_m}\cap \ga \;\oi{\alpha^\sigma}\alpha\,\big)\;.
$$ By Lemma \ref{lem:coincidencem}, for every $n\in\NN$, there exists
$\beta_n\in\Theta_\alpha$ such that the geodesic line
$\oi{\beta_n^\sigma}{\beta_n}$ coincides with the geodesic line
$\oi{\infty}{f_m}$ inside the horoballs $H_{n(m+1)+i+1}(f_m)$ for $1\leq
i\leq m$. Therefore $n(f_m,\alpha)\geq 2m$. For every such $\beta_n$,
the height of the geodesic line $\oi{\beta_n^\sigma}{\beta_n}$ inside
$H_{n(m+1)+i+1}(f_m)$ for $i=0,m+1$ is at most $1$. Hence for $n$ large
enough, the geodesic lines $\oi{\infty}{f_m}$ and
$\oi{\beta_n^\sigma}{\beta_n}$ coincide at most in one edge inside
$H_{n(m+1)+i+1}(f_m)$ for $i=0,m+1$. Thus $n(f_m,\alpha)\leq 2m+2$. Since
the stabiliser in $\Ga$ of both $*$ and the edge $e_\infty$ of $\lk(*)$
pointing towards $\infty$ acts transitively on the edges of $\lk(*)$
different from $e_\infty$, it is possible to adjust $\beta_n$ so that
$\oi{\infty}{f_m}$ and $\oi{\beta_n^\sigma}{\beta_n}$ do coincide in
exactly an edge inside $H_{n(m+1)+i+1}(f_m)$ for both $i=0,m+1$. Thus
$n(f_m,\alpha)= 2m+2$, as wanted.

Similarly, for every $n\in\NN$, consider the elements
$\beta'_n\in\Theta_\alpha$ such that the geodesic line
$\oi{{\beta'_n}^\sigma}{\beta'_n}$ coincides with the geodesic line
$\oi{\infty}{f'_m}$ inside the horoballs $H_{n(m+2)+i+1}(f'_m)$ for
$1\leq i\leq m$. Using the fact that the upper unipotent action of
$\FF_q[X]$ is simply transitive on the edges going out of the horoball
$\H_\infty$, and that two polynomials have the same action on
$\partial\H_\infty$ if and only if they differ by a constant, it is
possible to adjust $\beta'_n$ so that for $n$ large enough
$\oi{\infty}{f'_m}$ and $\oi{{\beta'_n}^\sigma}{\beta'_n}$ do coincide
in exactly one edge inside $H_{n(m+2)+1}(f'_m)$ and two edges inside
$H_{n(m+2)+m+1}(f'_m)$, but do not coincide in an edge inside
$H_{n(m+2)+m+2}(f'_m)$.  Thus $n(f'_m,\alpha)= 2m+3$, as wanted.
\cqfd

\medskip 
The next result shows that there are examples of quadratic
Lagrange spectra which contain a gap, that is, are not always of the
form $\{0\}\cup\big\{q^{-n}:n\in\NN, n\geq N\big\}$ for some
$N\in\NN$

\bprop\label{prop:gaps} Let $k\in\NN-\{0,1\}$, let $\{b_1,\dots,
b_{N'}\}$ be the set\footnote{with an arbitrary order and
  $N'=q^{2k+1}-q^k$} of the polynomials of degree at least $k$ and at
most $2k$, and let $\F_{k}=\bigcup_{\sigma\in{\mathcal S}_{N'}}
\Ga\cdot [\;\overline{b_{\sigma(1)}, \dots, b_{\sigma(N')}}\;]$. For
every $\alpha\in \F_{k}$, the number $q^{-2k+1}$ does not belong to
$\Sp(\alpha)$.  
\eprop

\medskip\dem Let $f\in \wh K-(K\cup\Theta_\alpha)$. Let us prove that
$c_\alpha(f)\neq q^{-2k+1}$, which gives the result.  There are three
cases to consider. 

Assume first that there exists $i_0\in\NN$ such that $\deg a_i(f)< k$
for all $i\ge i_0$. By Lemma \ref{lem:M2} (2), we have $c_\alpha(f)\ge
q^{-M_2(f)}\ge q^{-2(k-1)}=q^{-2k+2}$, and in particular $c_\alpha(f)
\neq q^{-2k+1}$.

Assume then that there exists a subsequence of coefficients
$a_{i_n}(f)$ for $n\in\NN$ such that $i_0\geq 1$ and $k\leq\deg
a_{i_n}(f)\le 2k$. Then $a_{i_n}(f)\in\{b_1,\dots, b_{N'}\}$, and
again there exists an element $\beta_{i_n}\in \Theta_\alpha$ for which
the intersection $\oi{\beta_{i_n}^\sigma} {\beta_{i_n}} \cap\oi \infty
f$ has length at least $2\deg a_{i_n}(f)\ge 2k$. Hence $c_\alpha(f)\le
q^{-2k}$, and in particular $c_\alpha(f)\neq q^{-2k+1}$.

If neither of the previous two cases occurs, we have $\deg a_i(f)>2k$
for $i$ large enough. By Lemma \ref{lem:inspectrum}, we have
$c_\alpha(f)\le q^{-4k}$, and in particular $c_\alpha(f)\neq
q^{-2k+1}$.  
\cqfd

\bcoro\label{coro:gaps} 
Let $\alpha$ be as in Proposition \ref{prop:gaps}. Then $\Sp(\alpha)$
contains a gap.  
\ecoro

\dem By Corollary \ref{coro:hurwitzconst}, we have $\max\Sp(\alpha)
=q^{-2}$, and $\Sp(\alpha)$ contains $q^{-n}$ for all $n$ large enough
by Theorem \ref{theo:hallfunctionfield}.  Thus, the spectrum has a
gap that contains $q^{-2k+1}$ by Proposition \ref{prop:gaps}.  \cqfd

{\small \bibliography{../biblio} }

\begin{thebibliography}{BAPP}

\bibitem[BaL]{BasLub01}
H.~Bass and A.~Lubotzky.
\newblock {\it Tree lattices}.
\newblock {Prog. in Math. {\bf 176}, Birkhäuser, 2001}.

\bibitem[BeN]{BerNak00}
V.~Berthé and H.~Nakada.
\newblock {\it On continued fraction expansions in positive characteristic:
  equivalence relations and some metric properties}.
\newblock {Expo. Math. {\bf 18} (2000) 257--284}.

\bibitem[BH]{BriHae99}
M.~R. Bridson and A.~Haefliger.
\newblock {\it Metric spaces of non-positive curva\-tu\-re}.
\newblock {Grund. math. Wiss. {\bf 319}, Springer Verlag, 1999}.

\bibitem[BPP]{BroParPau18}
A.~Broise-Alamichel, J.~Parkkonen, and F.~Paulin.
\newblock {\it Equidistribution and counting under equilibrium states in
   in negative curvature and trees. Applications to
  non-Archimedean Diophantine approximation}.
\newblock {Book preprint (365 pages) {\tt arXiv:1612.06717}, with an
Appendix by J.~Buzzi, to appear in Progress in Mathematics}.

\bibitem[Bug1]{Bugeaud14}
Y.~Bugeaud.
\newblock {\it On the quadratic Lagrange spectrum}.
\newblock {Math. Z. {\bf 276} (2014) 985--999}.

\bibitem[Bug2]{Bugeaud18}
Y.~Bugeaud.
\newblock {\it Nonarchimedean quadratic Lagrange spectra and 
continued fractions in power series fields}.
\newblock {Preprint {\tt arXiv:1804.03566}}.

\bibitem[CF]{CusFla89}
T.~Cusick and M.~Flahive.
\newblock {\it The Markoff and Lagrange spectra}.
\newblock {Math. Surv. Mono. {\bf 30}, Amer. Math. Soc. 1989}.

\bibitem[Gos]{Goss96}
D.~Goss.
\newblock {\it Basic structures of function field arithmetic}.
\newblock {Erg. Math. Grenz. {\bf 35}, Springer Verlag 1996}.

\bibitem[HP]{HerPau10}
S.~Hersonsky and F.~Paulin.
\newblock {\it On the almost sure spiraling of geodesics in negatively curved
  manifolds}.
\newblock {J. Diff. Geom. {\bf 85} (2010) 271--314}.

\bibitem[Khi]{Khinchin64}
A.~Y. Khinchin.
\newblock {\it Continued fractions}.
\newblock {Univ. Chicago Press, 1964}.

\bibitem[Las]{Lasjaunias00}
A.~Lasjaunias.
\newblock {\it A survey of Diophantine approximation in fields of power
  series}.
\newblock {Monat. Math. {\bf 130} (2000) 211--229}.

\bibitem[Lin]{Lin18}
X.~Lin.
\newblock {\it Quadratic Lagrange spectrum: I}.
\newblock {Math. Z. {\bf 289} (2018) 515--533}.

\bibitem[PaP1]{ParPau07}
J.~Parkkonen and F.~Paulin.
\newblock {\it Sur les rayons de Hall en approximation diophantienne}.
\newblock {Comptes Rendus Math. {\bf 344} (2007) 611--614}.

\bibitem[PaP2]{ParPau09}
J.~Parkkonen and F.~Paulin.
\newblock {\it On the closedness of approximation spectra}.
\newblock {J. Th. Nb. Bordeaux {\bf 21} (2009) 701--710}.

\bibitem[PaP3]{ParPau10GT}
J.~Parkkonen and F.~Paulin.
\newblock {\it Prescribing the behaviour of geodesics in negative curvature}.
\newblock {Geom. \& Topo. {\bf 14} (2010) 277--392}.

\bibitem[PaP4]{ParPau11MZ}
J.~Parkkonen and F.~Paulin.
\newblock {\it Spiraling spectra of geodesic lines in negatively curved
  manifolds}.
\newblock {Math. Z. {\bf 268} (2011) 101--142, Erratum: Math. Z. {\bf 276}
  (2014) 1215--1216}.

\bibitem[Pau]{Paulin02}
F.~Paulin.
\newblock {\it Groupe modulaire, fractions continues et approximation
  diophantienne en caract\'eristique $p$}.
\newblock {Geom. Dedi. {\bf 95} (2002) 65--85}.

\bibitem[Pej]{Pejkovic16}
T.~Pejkovi\'{c}.
\newblock {\it Quadratic Lagrange spectrum}.
\newblock {Math. Z. {\bf 283} (2016) 861--869}.

\bibitem[Poi]{Poitou53}
G.~Poitou.
\newblock {\it Sur l'approximation des nombres complexes par les nombres des
  corps imaginaires quadratiques dénués d'idéaux non principaux
  particulièrement lorsque vaut l'algorithme d'Euclide}.
\newblock {Ann. Scien. Ec. Norm. Sup. {\bf 70} (1953) 199--265}.

\bibitem[Ros]{Rosen02}
M.~Rosen.
\newblock {\it Number theory in function fields}.
\newblock {Grad. Texts Math. {\bf 210}, Springer-Verlag, 2002}.

\bibitem[Sch1]{Schmidt69}
A.~L. Schmidt.
\newblock {\it Farey simplices in the space of quaternions}.
\newblock {Math. Scand. {\bf 24} (1969) 31--65}.

\bibitem[Sch2]{Schmidt75}
A.~L. Schmidt.
\newblock {\it Diophantine approximation of complex numbers}.
\newblock {Acta Math. {\bf 134} (1975) 1--85}.

\bibitem[Sch]{Schmidt00}
W.~Schmidt.
\newblock {\it On continued fractions and diophantine approximation in power
  series fields}.
\newblock {Acta Arith.~{\bf XCV} (2000) 139--166}.

\bibitem[Ser]{Serre83}
J.-P. Serre.
\newblock {\it Arbres, amalgames, SL$_2$}.
\newblock {3\`eme \'ed. corr., Ast\'erisque {\bf 46}, Soc. Math. France, 1983}.

\end{thebibliography}

\bigskip
{\small
\noindent \begin{tabular}{l} 
Department of Mathematics and Statistics, P.O. Box 35\\ 
40014 University of Jyv\"askyl\"a, FINLAND.\\
{\it e-mail: jouni.t.parkkonen@jyu.fi}
\end{tabular}
\medskip

\noindent \begin{tabular}{l}
Laboratoire de mathématique d'Orsay,\\
UMR 8628 Univ. Paris-Sud et CNRS,\\
Universit\'e Paris-Saclay,\\
91405 ORSAY Cedex, FRANCE\\
{\it e-mail: frederic.paulin@math.u-psud.fr}
\end{tabular}
}

\end{document}